\documentclass[12pt]{article}
\usepackage{amssymb}
\usepackage{amsmath}
\usepackage{amsopn}
\numberwithin{equation}{section}

\newtheorem{theorem}{Theorem}
\newtheorem{proposition}{Proposition}[section] 
\newtheorem{corollary}[proposition]{Corollary}
\newtheorem{lemma}[proposition]{Lemma}

\let\ds=\displaystyle

\def\bn{\begin{equation}}
\def\ed{\end{equation}}

\def\f{F^+}
\def\e{E^+}

\newcommand{\eain}[2]{e_{#1}[#2]}
\newcommand{\fain}[2]{f_{#1}[#2]}

\def\<{\langle}
\def\>{\rangle}
\newcommand{\CC}{{\mathbb C}}

\newcommand{\ZZ}{{\mathbb Z}}

\let\z=\ZZ
\def\r#1{(\ref{#1})}
\let\rf=\r
\def\ot{\otimes}
\def\a{\alpha}   
\def\b{\beta} 

\def\al{\alpha}

\newcommand{\Uqdva}{{U_{q}^{}(\widehat{\mathfrak{sl}}_{2})}}
\def\sk#1{\left(#1\right)}

\def\i{\iota}

\def\ggg{\mathfrak{g}}
\def\slg{\mathfrak{sl}}

\def\Uqgln{U_q(\widehat{\mathfrak{gl}}_N)}
\def\Uqgl#1{U_q(\widehat{\mathfrak{gl}}_{#1})}

\def\Uqsln{U_q(\widehat{\mathfrak{sl}}_N)}
\def\Uqsl2{U_q(\widehat{\mathfrak{sl}}_2)}

\def\ee{E}\def\ff{F}

\def\w{\mathbf{w}}
\def\prodl{\mathop{\overleftarrow\prod}\limits}

\def\ct{\mathbb{T}}
\def\End{\textrm{End}}
\def\Uqbp{U_q(\mathfrak{b}_+)}
\def\bbb{\mathbb{B}}

\def\ccR{\mathbb{R}}

\def\ord{\prec}
\def\E{{\sf e}}

\begin{document}

\begin{center}

\hfill ITEP-TH-50/06\\
\hfill math.QA/0610517\\
\bigskip
{\Large\bf Off-shell Bethe vectors and Drinfeld currents}
\par\bigskip\medskip
{\bf
S.~Khoroshkin$^{\star}$\footnote{E-mail: khor@itep.ru},\ \
S.~Pakuliak$^{\star\bullet}$\footnote{E-mail: pakuliak@theor.jinr.ru},\ \
V.~Tarasov$^{*\diamond}$\footnote{E-mail: vt@pdmi.ras.ru, vtarasov@math.iupui.edu}}
\par\bigskip\medskip
$^\star${\it Institute of Theoretical \& Experimental Physics, 117259
Moscow, Russia}
\par\smallskip
$^\bullet${\it Laboratory of Theoretical Physics, JINR,
141980 Dubna, Moscow reg., Russia}
\par\smallskip
$^*${\it St.Petersburg Branch of Steklov Mathematical Institute\\
Fontanka 27, St.Peterburg 191011, Russia}
\par\smallskip
$^\diamond${\it Department of Mathematical Sciences, IUPUI, Indianapolis, IN
46202, USA}
\bigskip

\end{center}

\thispagestyle{empty}

\begin{abstract}
In this paper we compare two constructions of weight functions (off-shell Bethe
vectors) for the quantum affine algebra $U_q(\widehat{\mathfrak{gl}}_N)$.
The first construction comes from the algebraic nested Bethe ansatz. The second
one is defined in terms of certain projections of products of Drinfeld
currents. We show that two constructions give the same result in tensor
products of vector representations of $U_q(\widehat{\mathfrak{gl}}_N)$.
\end{abstract}


\section{Introduction}

Off-shell Bethe vectors in integrable models associated with the Lie algebra
$\mathfrak{gl}_N$ have appeared in \cite{KR83} in the framework of the algebraic
nested Bethe ansatz. For $N=2$ they have the form $B(t_1)\cdots B(t_k)v$, where
$B(u)=T_{12}(u)$ is an element of the monodromy matrix and $v$ is the highest
weight vector of an irreducible finite-dimensional representation of $\Uqgl2$.
For $N>2$, off-shell Bethe vectors are defined in \cite{KR83} inductively.
They are functions of several complex variables $t^1_1\dots t^{N-1}_k$ labeled
by two indices, the superscript corresponding to a simple root of
${\mathfrak{sl}}_N$. If the variables $t^1_1\dots t^{N-1}_k$ satisfy the Bethe
ansatz equations, the Bethe vectors are eigenvectors of the transfer matrix of
the system.

Off-shell Bethe vectors also serve for integral representations of solutions to
the $q$-difference Knizhnik-Zamolodchikov (qKZ) equations \cite{R, VT,VT1}.
In this case they are known under the name of ``weight function''. It has been
observed in \cite{VT} that weight functions have very special comultiplication
properties that allow one to express a weight function in a tensor product of
representations in terms of weight functions in the tensor factors. The
comultiplication properties of weight functions are essential for constructing
solutions to the qKZ equations. In this paper we start from these properties
and define a weight function as a collection of rational functions with values
in representations of the quantum affine algebra satisfying suitable
comultiplication relations.

A new approach for construction of weight functions has been proposed recently
in \cite{KP,EKhP}. It is based on the ``new realization'' of quantum affine
algebras \cite{D88}. In this approach, the key role is played by certain
projections to the intersection of Borel subalgebras of different kind of
the quantum affine algebra. Those projections were introduced in \cite{ER} and
were used in \cite{DKhP} to obtain integral formulae for the universal
$\mathcal{R}$-matrix of the quantum affine algebra $\Uqdva$. It is shown in
\cite{KP,EKhP} that acting by a projection of a product of Drinfeld currents
on highest weight vectors of irreducible finite-dimensional representations
of $\Uqgln$ one obtains a collection of rational functions with the required
comultiplication properties, that is, a weight function.

In this paper we compare two constructions of weight functions for the quantum
algebra $\Uqgln$. We conjecture that the constructions give the same result
for any irreducible finite-dimensional representation of $\Uqgln$. We prove
this conjecture for tensor products of the vector representations of $\Uqgln$.
To this end, we show that weight functions defined by the projections and
those given by the algebraic Bethe ansatz satisfy the same recurrence relations
with respect to the rank $N$ of the algebra. To obtain the recurrence relations
we use a generalization of the Ding-Frenkel isomorphisms of two realizations of
$\Uqgln$.

The paper is organized as follows. In Section 2 we recall two descriptions of
$\Uqgln$: in terms of the fundamental $L$-operators and in terms of Drinfeld
currents. We define weight functions and symmetric (or modified) weight
functions by their coalgebraic properties. We pay special attention to the
symmetry properties of weight functions. In Section 3 we describe the
construction of a weight function by means of projections of Drinfeld currents.
This is done by applying the construction of \cite{EKhP} to the subalgebra
$\Uqsln$ in $\Uqgln$. In Section 4 we describe, following \cite{VT},
the $L$-operator construction of the $\Uqgln$ weight function. In Section 5
we prove the main result of the paper that the two constructions of weight
functions give the same result for tensor products of the vector
representations of $\Uqgln$. Besides, Section 5 contains a description of
projections of composed currents, which generalizes the Ding-Frenkel
isomorphism.

\section{Quantum affine algebra $\Uqgln$}

\subsection{$L$-operator description}

\label{section2.1}

Let $\E_{ij}\in{\textrm{End}}(\CC^N)$ be a matrix with the only nonzero entry
equal to $1$
at the intersection of the $i$-th row and $j$-th column.
Let $R(u,v)\in{\textrm{End}}(\CC^N\ot\CC^N)\ot \CC[[{v}/{u}]]$,
\begin{equation}\label{UqglN-R}
\begin{split}
R(u,v)\ =\ &\frac{qu-q^{-1}v}{u-v}\ \sum_{1\leq i\leq N}\E_{ii}\ot \E_{ii}\ +\
\sum_{1\leq i<j\leq N}(\E_{ii}\ot \E_{jj}+\E_{jj}\ot \E_{ii})
\\
+\ &\frac{q-q^{-1}}{u-v}\sum_{1\leq i<j\leq N}
(v \E_{ij}\ot \E_{ji}+ u \E_{ji}\ot \E_{ij})\,,
\end{split}
\end{equation}
be the standard trigonometric $R$-matrix associated with the vector
representation of ${\mathfrak{gl}}_N$. It satisfies the Yang-Baxter equation
\begin{equation}\label{YBeq}
R^{12}(u_1,u_2)R^{13}(u_1,u_3)R^{23}(u_2,u_3)=
R^{23}(u_2,u_3)R^{13}(u_1,u_3)R^{12}(u_1,u_2)\,,
\end{equation}
and the inversion relation
\begin{equation}\label{unitar}
R^{(12)}(u_1,u_2)R^{(21)}(u_2,u_1)=\frac{(qu_1-q^{-1}u_2)(q^{-1}u_1-qu_2)}{(u_1-u_2)^2}.
\end{equation}

The algebra $\Uqgln$ (with the zero central charge and the gradation operator
dropped out) is a unital associative algebra generated by the modes
$L^\pm_{ij}[\pm k]$, $k\geq 0$, $1\leq i,j\leq N$, of the $L$-operators
$L^\pm(z)=\sum_{k=0}^\infty\sum_{i,j=1}^N \E_{ij}\otimes L^\pm_{ij}[\pm k]z^{\mp
k}$, subject to relations
\begin{equation}\label{L-op-com}
\begin{split}
&R(u,v)\cdot (L^\pm(u)\ot \mathbf{1})\cdot (\mathbf{1}\ot L^\pm(v))=
(\mathbf{1}\ot L^\pm(v))\cdot (L^\pm(u)\ot \mathbf{1})\cdot R(u,v)\\
&R(u,v)\cdot (L^+(u)\ot \mathbf{1})\cdot (\mathbf{1}\ot L^-(v))=
(\mathbf{1}\ot L^-(v))\cdot (L^+(u)\ot \mathbf{1})\cdot R(u,v),\\
&L^+_{ij}[0]=L^-_{ji}[0]=0,\qquad L^+_{kk}[0]L^-_{kk}[0]=1,
\qquad 1\leq i<j \leq N,
\quad 1\leq k\leq N\ .
\end{split}
\end{equation}

The coalgebraic structure of the algebra $\Uqgln$ is defined by the rule
\begin{equation}\label{coprL}
\Delta \sk{L^\pm_{ij}(u)}=\sum_{k=1}^N\ L^\pm_{kj}(u)\otimes
L^\pm_{ik}(u)\,.
\end{equation}

\subsection{The current realization of $\Uqgln$}

The algebra $\Uqgln$ in the current realization (with the zero central charge
and the gradation operator dropped out) is generated by the modes of the Cartan
currents
$$k_i^\pm(z)=\sum_{m\geq0}k_i^\pm[\pm m]z^{\mp m},\qquad k^+_i[0]k^-_i[0]=1\,,$$
$i=1,\ldots,N$, and by the modes of the generating functions 
\begin{equation}\label{currents}
\ee_i(z)=\sum_{n\in\ZZ}\ee_i[n]z^{-n,}\qquad
\ff_i(z)=\sum_{n\in\ZZ}\ff_i[n]z^{-n}\, ,
\end{equation}
$i=1,\ldots,N-1$, subject to relations
$$
(q^{-1}z-q^{}w)\ee_{i}(z)\ee_{i}(w)=
\ee_{i}(w)\ee_{i}(z)(q^{}z-q^{-1}w)\, ,
$$
$$
(z-w)\ee_{i}(z)\ee_{i+1}(w)=
\ee_{i+1}(w)\ee_{i}(z)(q^{-1}z-qw)\, ,
$$
$$
(q^{}z-q^{-1}w)\ff_{i}(z)\ff_{i}(w)=
\ff_{i}(w)\ff_{i}(z)(q^{-1}z-q^{}w)\, ,
$$
$$
(q^{-1}z-qw)\ff_{i}(z)\ff_{i+1}(w)=
\ff_{i+1}(w)\ff_{i}(z)(z-w)\, ,
$$
$$
k_i^\pm(z)\ff_i(w)\left(k_i^\pm(z)\right)^{-1}=
\frac{q^{-1}z-qw}{z-w}\ff_i(w)\, ,
$$
\begin{equation}\label{gln-com}
k_{i+1}^\pm(z)\ff_i(w)\left(k_{i+1}^\pm(z)\right)^{-1}=
\frac{q^{}z-q^{-1}w}{z-w}\ff_i(w)\, ,
\end{equation}
$$
k_i^\pm(z)\ff_j(w)\left(k_i^\pm(z)\right)^{-1}=\ff_j(w)
\qquad {\rm if}\quad i\not=j,j+1\, ,
$$
$$
k_i^\pm(z)\ee_i(w)\left(k_i^\pm(z)\right)^{-1}=
\frac{z-w}{q^{-1}z-q^{}w}\ee_i(w)\, ,
$$
$$
k_{i+1}^\pm(z)\ee_i(w)\left(k_{i+1}^\pm(z)\right)^{-1}=
\frac{z-w}{q^{}z-q^{-1}w}\ee_i(w)\, ,
$$
$$
k_i^\pm(z)\ee_j(w)\left(k_i^\pm(z)\right)^{-1}=\ee_j(w)
\qquad {\rm if}\quad i\not=j,j+1\, ,
$$
$$
[\ee_{i}(z),\ff_{j}(w)]=
\delta_{{i},{j}}\ \delta(z/w)\ (q-q^{-1})\left(
k^+_{i}(z)/k^+_{i+1}(z)-k^-_{i}(w)/k^-_{i+1}(w)\right)\, ,
$$
together with the Serre relations
\begin{equation}
\begin{split}
{\rm Sym}_{z_1,z_{2}}
(\ee_{i}(z_1)\ee_{i}(z_2)\ee_{i\pm 1}(w)
&-(q+q^{-1})\ee_{i}(z_1)\ee_{i\pm 1}(w)\ee_{i}(z_2)+\\
&+\ee_{i\pm 1}(w)\ee_{i}(z_1)\ee_{i}(z_2))=0\, ,\\
\label{serre}
{\rm Sym}_{z_1,z_{2}}
(\ff_{i}(z_1)\ff_{i}(z_2)\ff_{i\pm 1}(w)
&-(q+q^{-1})\ff_{i}(z_1)\ff_{i\pm 1}(w)\ff_{i}(z_2)+\\
&+\ff_{i\pm 1}(w)\ff_{i}(z_1)\ff_{i}(z_2))=0\, .
\end{split}
\end{equation}

To construct an isomorphism between the $L$-operator and current realizations
of the algebra $\Uqgln$, one has to decompose the $L$-operators into
the Gauss coordinates
\begin{equation}\label{L-op}
L^\pm(z)=\sk{\sum_{i=1}^N \E_{ii}+\sum^N_{i<j}F^\pm_{i,j}(z)\E_{ij}}
\cdot\sk{\sum^N_{i=1}k^\pm_{i}(z)\E_{ii}}\cdot
\sk{\sum_{i=1}^N \E_{ii}+\sum^N_{i<j}E^\pm_{j,i}(z)\E_{ji}}
\end{equation}
and for $i=1,\ldots,N-1$ to identify the total currents and the linear combinations
of the nearest to the diagonal Gauss coordinates \cite{DF}
\begin{equation}\label{DF-iso}
E_i(z)=\ee^+_{i+1,i}(z)-\ee^-_{i+1,i}(z)\,,\quad
F_i(z)=\ff^+_{i,i+1}(z)-\ff^-_{i,i+1}(z)\ .
\end{equation}

The diagonal Gauss coordinates of the $L$-operators coincide with the Cartan
currents $k^\pm_i(z)$ and are denoted by the same letter. The results of
\cite{DF} say nothing about relations between Gauss coordinates
$\ff^\pm_{i,j}(z)$ and $\ee^\pm_{j,i}(z)$ of the $L$-operators for $j-i>1$ and
the currents $F_i(z)$, $E_i(z)$. Some of these relations, that we will need for
our construction, will be described in subsection~\ref{MoreDF}.

In \cite{D88} the \emph{current} Hopf structure for the algebra $\Uqgln$ has
been defined,
\begin{equation}\label{gln-copr}
\begin{split}
\Delta^{(D)}\sk{E_i(z)}&=E_i(z)\ot 1 + k^-_{i}(z)\sk{k^-_{i+1}(z)}^{-1}\ot
E_i(z),\\
\Delta^{(D)}\sk{F_i(z)}&=1\ot F_i(z) + F_i(z)\ot
k^+_{i}(z)\sk{k^+_{i+1}(z)}^{-1},\\
\Delta^{(D)}\sk{k^\pm_i(z)}&=k^\pm_i(z)\ot k^\pm_{i}(z).
\end{split}
\end{equation}

We consider two types of Borel subalgebras of the
algebra $\Uqgln$. Borel subalgebras $U_q(\mathfrak{b}_\pm)\subset \Uqgln $
are generated by the modes of the $L$-operators $L^\pm(z)$, respectively.

Another type of Borel subalgebras is related to the current realization of
$\Uqgln$. The Borel subalgebra $U_F\subset \Uqgln$ is generated by
modes of the currents
$F_i[n]$, $k^+_j[m]$, $i=1,\ldots,N-1$, $j=1,\ldots,N$, $n\in\ZZ$
and $m\geq0$. The Borel subalgebra
$U_E\subset \Uqgln$ is generated by modes of the currents
$E_i[n]$, $k^-_j[-m]$, $i=1,\ldots,N-1$, $j=1,\ldots,N$, $n\in\ZZ$ and
$m\geq0$. We will consider also a subalgebra $U'_F\subset U_F$, generated
by the elements
$F_i[n]$, $k^+_j[m]$, $i=1,\ldots,N-1$, $j=1,\ldots,N$, $n\in\ZZ$
and $m>0$, and a subalgebra $U'_E\subset U_E$ generated by
the elements
$E_i[n]$, $k^-_j[-m]$, $i=1,\ldots,N-1$, $j=1,\ldots,N$, $n\in\ZZ$
and $m>0$.
Further, we will be interested in the intersections,
\begin{equation}
\label{Intergl}
U_f^-=U'_F\cap U_q(\mathfrak{b}_-)\,,\qquad
U_F^+=U_F\cap U_q(\mathfrak{b}_+)\,
\end{equation}
and will describe properties of projections to these intersections.

\subsection{ A weight function}\label{section2.3}

We call a vector $v$ \emph{a weight singular vector} if it is annihilated
by any mode of the currents $E_i[n]$, $i=1,\ldots,N-1$, $n\geq 0$ and is
an eigenvector for the Cartan currents $k^+_i(z)$, $i=1,\ldots,N$
\begin{equation}\label{hwv}
E^+_{i+1,i}(z)\cdot v=0\ ,\qquad k^+_i(z)\cdot v=\Lambda_i(z)\,v\,,
\end{equation}
where $\Lambda_i(z)$ is a meromorphic function, decomposed as a power series
in $z^{-1}$.
The $L$-operator \r{L-op}, acting on a weight singular vector $v$, become
upper-triangular
\begin{equation}\label{L-op-tr}
L^+_{ij}(z)\ v=0\,,\quad i>j\,,\quad L^+_{ii}(z)\ v=\Lambda_i(z)\ v\,,\quad
i=1,\ldots,N\,.
\end{equation}

We define a weight function by its comultiplication properties.

Let $\Pi$ be the set $\{1,\ldots,N-1\}$
of indices of simple positive roots of $\mathfrak{gl}_N$.
A finite collection $I=\{i_1,\dots,i_n\}$ with a linear ordering
$i_i\ord\cdots\ord i_n$ and a map $\iota:I\to\Pi$ is called an {\it ordered
$\Pi$-multiset}. Sometimes, we denote the map $\iota$ by $\iota_I$ and call it
a ``{\it colouring map}''. A morphism between two ordered $\Pi$-multisets $I$
and $J$ is a map $m:I\to J$ that respects the orderings in $I$ and $J$
and intertwines the colouring maps: $\iota_J m=m\iota_I$. In particular, any
subset $I'\subset I$ of a $\Pi$-ordered multiset has a unique structure of
$\Pi$-ordered multiset, such that the inclusion map is a morphism of
$\Pi$-ordered multisets.

To each $\Pi$-ordered multiset $I=\{i_1,\dots,i_n\}$ we attach an ordered set
of variables $\{t_i|i\in I\}=\{t_{i_1},\dots,t_{i_n}\}$. Each variable has its
own ``colour'': $\iota(i_k)\in\Pi$.

Let $i$ and $j$ be elements of some ordered $\Pi$-multiset. Define a
rational function
\begin{equation}\label{gam}
\gamma(t_i,t_j)=\left\{
\begin{array}{ll}
\ds\frac{qt_i-q^{-1}t_j}{t_i-t_j}\ ,\quad&\mbox{if}\quad \i(i)=\i(j)+1\ ,\\[5mm]
\ds\frac{t_i-t_j}{q^{-1}t_i-qt_j}\ ,\quad&\mbox{if}\quad \i(j)=\i(i)+1\ ,\\[5mm]
\ds\frac{q^{-1}t_i-qt_j}{qt_i-q^{-1}t_j}\ ,\quad&\mbox{if}\quad \i(i)=\i(j)\ ,\\
[5mm]
1\ ,\quad&\mbox{otherwise}\ .
\end{array}
\right.
\end{equation}

Assume that for any representation $V$ of $\Uqgln$, generated by a weight
singular vector $v$, and any ordered $\Pi$-multiset $I=\{i_1,\dots i_n\}$,
there is a $V$-valued rational function $w_{V,I}(t_{i_1},\ldots,t_{i_n})\in V$
depending on the variables $\{t_i|i\in I\}$. We call such a collection of
rational functions a {\it weight function} $w$, if:
\begin{itemize}
\item[(a)] The rational function, corresponding to the empty set,
is equal to $v$,
\begin{equation}\label{empty}
w_{V,\emptyset}\equiv v\ .
\end{equation}
\item[(b)] The function $w_{V,I}(t_{i_1},\ldots,t_{i_n})$ depends only on
an isomorphism class of an ordered $\Pi$-multiset, that is, for any isomorphism
$f:I\to J$ of ordered $\Pi$-multisets we have
\begin{equation}
\label{isomorph}
w_{V,I}(t_{f(i)}|_{i\in I})=w_{V,J}(t_j|_{j\in J})\,.
\end{equation}
\item[(c)]
The functions $w_{V,I}$ satisfy the following comultiplication property.
Let $V=V_1\otimes V_2$ be a tensor product of two
representations generated by the singular vectors $v_1$, $v_2$ and
weight series $\{\Lambda_b^{(1)}(u)\}$ and
$\{\Lambda_b^{(2)}(u)\}$, $b=1,\ldots,N$. Then for any multiset $I$ we
have
\end{itemize}
\begin{equation}\label{weight2}
w_{V,I}(t_{i} |_{i\in I})=
\sum\limits_{I=I_1\coprod I_2}
w_{V_1,I_1}(t_{i}|_{i\in I_1})\otimes w_{V_2,I_2}(t_{i}|_{i\in I_2})
\ \Phi_{I_1,I_2}(t_{i}|_{i\in I})
\prod\limits_{{j\in I_1}}
\frac{\Lambda^{(2)}_{\iota(j)}(t_{j})}{\Lambda^{(2)}_{\iota(j)+1}(t_{j})}
\,,
\end{equation}
{where}
\begin{equation}
\Phi_{I_1,I_2}(t_{i}|_{i\in I})=
\prod\limits_{\substack{i\in I_1,\ j\in I_2\\ i\ord j}}
\gamma(t_i,t_j).
\end{equation}
The summation in \r{weight2} runs over all possible decompositions of
the ordered multiset $I$ into a disjoint union of two non-intersecting ordered
submultisets $I_1$ and $I_2$.

Note that the comultiplication property relation \r{weight2} is not a recurrence
relation, that is, it does not allow us to reconstruct functions $w_{V,I}$
for all ordered multisets $I$ starting from the functions which correspond to
the multisets with $|I|=1$.
\bigskip

Let $I=\{i_1,\dots,i_n\}$ and $J=\{j_1,\dots,j_n\}$ be two ordered
$\Pi$-multisets. Let $\sigma:I\to J$ be an invertible map, which intertwines
the colouring maps, $\iota_J \sigma=\sigma\iota_I$, but does not necessarily
respect the orderings in $I$ and $J$ (that is, $\sigma$ is a ``permutation'' on
classes of isomorphisms of ordered multisets).

Let $w(t_j|_{j\in J})$ be a function of the variables $t_j|_{j\in J}$.
Define a {\it pullback} $^{\sigma,\gamma}\!w(t_i|_{i\in I})$ by the rule
\begin{equation}\label{pullback}
^{\sigma,\gamma}\!w(t_i|_{i\in I})=w(t_{\sigma(i)}|_{i\in I})
\prod_{\substack{i,j\in I\\ i\ord j,\
\sigma(j)\ord \sigma(i)}}\gamma(t_i,t_j)\,.
\end{equation}

Let $I$ and $I'$ be ordered $\Pi$-multisets, and $\sigma:I\to I'$ an invertible
map, intertwining the colouring maps. Then its restriction to any subset
$J\subset I$ is an invertible map of $J$ to $\sigma(J)$, intertwining the
colouring maps.
\begin{proposition} \label{proppull}
Let $w$ be a weight function, $I,J$ ordered $\Pi$-multisets, and
$\sigma:I\to J$ an invertible map, intertwining the colouring maps.
Then we have
\begin{equation*}
\begin{split}
^{\sigma,\gamma}\, & \!w_{V,J}(t_{i} |_{i\in I})={}
\\[4pt]
& {}=\sum\limits_{I=I_1\coprod I_2}
{}^{\sigma,\gamma}\!w_{V_1,\sigma(I_1)}(t_{i}|_{i\in I_1})
\otimes\
^{\sigma,\gamma}\!w_{V_2,\sigma(I_2)}(t_{i}|_{i\in I_2})
\ \Phi_{I_1,I_2}(t_{i}|_{i\in I})
\prod\limits_{{j\in I_1}}
\frac{\Lambda^{(2)}_{\iota(j)}(t_{j})}{\Lambda^{(2)}_{\iota(j)+1}(t_{j})}
\,.
\end{split}
\end{equation*}
\end{proposition}
\noindent
The proposition means that the pullback operation \rf{pullback} is compatible
with the comultiplication rule \rf{weight2}.

We call a weight function $w$ {\it $q$-symmetric}, if for any ordered
$\Pi$-multisets $I$ and $J$ and an invertible map $\sigma:I\to J$,
intertwining the colouring maps, we have
$$^{\sigma,\gamma}w_{V,J}(t_{i} |_{i\in I})= w_{V,I}(t_{i} |_{i\in I})\,.$$

\subsection{A modified weight function}

Given elements $i,j$ of some ordered multiset define two functions
$\tilde{\gamma}(t_i,t_j)$ and $\beta(t_i,t_j)$ by the formulae
\begin{equation*}
\tilde\gamma(t_i,t_j)=\left\{
\begin{array}{ll}
\ds\frac{t_i-t_j}{qt_i-q^{-1}t_j}\ ,\quad&\mbox{if}\quad \i(i)=\i(j)+1\ ,\\[5mm]
\ds\frac{q^{-1}t_i-qt_j}{t_i-t_j}\ ,\quad&\mbox{if}\quad \i(j)=\i(i)+1\ ,\\[5mm]
1\ ,\quad&\mbox{otherwise}
\end{array}
\right.
\end{equation*}
and
\begin{equation}\label{beta}
\beta(t_i,t_j)=\left\{
\begin{array}{ll}
\ds\frac{q^{-1}t_i-qt_j}{t_i-t_j}\ ,\quad&\mbox{if}\quad \i(i)=\i(j)\ ,\\[5mm]
1\ ,\quad&\mbox{otherwise\,.}
\end{array}
\right.
\end{equation}

A collection of rational $V$-valued functions $\w_{V,I}(t_{i}|_{i\in I})$,
depending on a representation $V$ of $\Uqgln$, generated by a weight singular
vector $v$, and an ordered $\Pi$-multiset $I$, is called a {\it modified weight
function} $\w$, if it satisfies conditions (a), (b), see \r{empty},
\r{isomorph}, and condition (c$'$):

\begin{itemize}
\item[(c$'$)]
Let $V=V_1\otimes V_2$ be a tensor product of two representations generated by
the singular vectors $v_1$, $v_2$ and weight series $\{\Lambda_b^{(1)}(u)\}$
and $\{\Lambda_b^{(2)}(u)\}$, $b=1,\ldots,N$. Then for any multiset $I$ we have
\end{itemize}
\begin{equation}\label{weight22}\begin{split}
\w_{V,I}(t_{i} |_{i\in I})
=&\sum\limits_{I=I_1\coprod I_2}
\w_{V_1,I_1}(t_{i}|_{i\in I_1})
\otimes
\w_{V_2,I_2}(t_{i}|_{i\in I_2})
\cdot\tilde{\Phi}_{I_1,I_2}(t_{i}|_{i\in I})\times\\
\times&\prod\limits_{{j\in I_1}}
{\Lambda^{(2)}_{\iota(j)}(t_{j})}
\prod\limits_{{j\in I_2}}
{\Lambda^{(1)}_{\iota(j)+1}(t_{j})},
\end{split}\end{equation}
where
\begin{eqnarray*}
\tilde{\Phi}_{I_1,I_2}(t_{i}|_{i\in I})=
\prod_{i\in I_1,\ j\in I_2} \beta(t_i,t_j)
\prod_{\substack{i\in I_2,\ j\in I_1\\ i\ord j}}\tilde\gamma(t_i,t_j)\,.
\end{eqnarray*}

Let $I=\{i_1,\dots,i_n\}$ and $J=\{j_1,\dots,j_n\}$ be two ordered
$\Pi$-multisets, and
$\sigma:I\to J$ an invertible map, intertwining the colouring maps,
$\iota_J \sigma=\sigma\iota_I$.
Let $\w(t_j|_{j\in J})$ be a function of the variables $t_j|_{j\in J}$.
Define a pullback $^{\sigma,\tilde{\gamma}}\w(t_i|_{i\in I})$ by the rule
\begin{equation}\label{pullbacka}
^{\sigma,\tilde{\gamma}}\w(t_i|_{i\in I})=\w(t_{\sigma(i)}|_{i\in I})
\prod_{\substack{i,j\in I\\ i\ord j,\ \sigma(j)\ord \sigma(i)}}
\tilde{\gamma}(t_i,t_j)\,.
\end{equation}

\begin{proposition} \label{proppulla}
Let $\w$ be a modified weight function, $I,J$ ordered $\Pi$-multisets,
and $\sigma:I\to J$ an invertible map, intertwining the colouring maps.
Then we have
\begin{equation*}\begin{split}
^{\sigma,\tilde{\gamma}}\w_{V,J}(t_{i} |_{i\in I})=
\sum\limits_{I=I_1\coprod I_2}{}
{}^{\sigma,\tilde{\gamma}}\w_{V_1,\sigma(I_1)}(t_{i}|_{i\in I_1})
\otimes\
^{\sigma,\tilde{\gamma}}\w_{V_2,\sigma(I_2)}(t_{i}|_{i\in I_2})\times{} &
\\[4pt]
{}\times\tilde{\Phi}_{I_1,I_2}(t_{i}|_{i\in I})
\prod\limits_{{j\in I_1}}{\Lambda^{(2)}_{\iota(j)}(t_{j})}
\prod\limits_{{j\in I_2}}
{\Lambda^{(1)}_{\iota(j)+1}(t_{j})}\, & .
\end{split}\end{equation*}
\end{proposition}
We call a modified weight function $\w$ {\it $q$-symmetric}, if for any two
ordered $\Pi$-multisets $I$ and $J$ and an invertible map $\sigma:I\to J$,
intertwining the colouring maps, we have
$$^{\sigma,\tilde{\gamma}}\w_{V,J}(t_{i} |_{i\in I})=
\w_{V,I}(t_{i} |_{i\in I})\,.$$

For an ordered $\Pi$-multiset $I=\{i_1,i_{2},\ldots,i_n\}$, let
$\bar I=\{i_n,i_{n-1},\ldots,i_1\}$ be the ordered $\Pi$-multiset with
the colouring map, $\i_{\bar I}(i_k)=\i_I(i_k)$, $k=1,\ldots,n$.
\begin{proposition}\label{prop2.3}${}$
\begin{itemize}
\item[(i)]
Let $w$ be a weight function. Then the collection
$\w_{V,I}(t_i|_{i\in I})$, where
\begin{equation}\label{swf}
\w_{V,I}(t_i|_{i\in I})=w_{V,\bar I}(t_i|_{i\in\bar I})
\prod\limits_{i\ord j } \beta(t_i,t_j)\,
\prod_{i\in I}
\Lambda_{\i(i)+1}(t_i)
\end{equation}
is a modified weight function.
\item[(ii)]
Let $\w$ be a modified weight function. Then the collection
$w_{V,I}(t_i|_{i\in I})$, where
\begin{equation*}
w_{V,I}(t_i|_{i\in I})=\w_{V,\bar I}(t_i|_{i\in\bar I})\,
\prod\limits_{i\ord j} \frac1{\beta(t_j,t_i)}\;
\prod_{i\in I}\frac1{\Lambda_{\i(i)+1}(t_i)}\;
\end{equation*}
is a weight function.
\item[(iii)]
If $w$ is a $q$-symmetric weight function, then $\w$ is $q$-symmetric
modified weight function, and vice versa.
\end{itemize}
\end{proposition}
The last proposition means that we have a bijection between weight functions
and modified weight functions.

\section{Weight functions and Drinfeld currents}

\subsection{Quantum affine algebra $\Uqsln$}
We are using two description of the quantum affine algebra $\Uqsln$:
in terms of Chevalley generators and the current realization.

The algebra $\Uqsln$ (with zero central charge and the grading element
dropped out) is generated by the Chevalley generators $e_{\pm\a_i}$,
$k^{\pm1}_{\a_i}$, where $i=0,1,\ldots,N-1$ and $\prod_{i=0}^N k_{\a_i}=1$,
subject to relations
\begin{equation}\label{Chevalley}
k_{\a_i}e_{\pm\al_j}k^{-1}_{i}\,=\,q_i^{\pm a_{ij}}e_{\pm\al_j}\,,
\quad
[e_{i},e_{-\al_j}]\,=\,
\delta_{ij}\frac{k_{\a_i}-k^{-1}_{\a_i}}{q_i-q_i^{-1}}\,,
\end{equation}
\begin{equation}\label{re2}
\begin{split}
&\sum_{r=0}^{m_{i,j}} {e^{(r)}_{\pm\alpha_i}} e_{\pm\alpha_j}
{e^{(m_{i,j}-r)}_{\pm\alpha_i}}=0,\quad\text{where}\quad
m_{i,j}=1-(\alpha_i,\alpha_j),\quad i\neq j,
\\
&e_{\pm\alpha_i}^{(r)}=\frac{e_{\pm\alpha_i}^r}{[k]_q!},\qquad
[k]_q!=[k]_q [k-1]_q\dots [2]_1 [1]_q,\qquad
[k]_q=\frac{q^k-q^{-k}}{q-q^{-1}},
\end{split}
\end{equation}
where
$a_{i,j}=(\al_i,\al_j)$ is the
Cartan matrix of the affine algebra $\widehat{\mathfrak{sl}}_N$,
$(\al_{i}, \al_{i})=2$, $(\al_{i}, \al_{j})=-1$, if $i-j=\pm 1\ \mod N$.

The comultiplication map
is given by the formulae:
\begin{equation}\label{copr}\begin{split}
\Delta(e_{\a_i})&=e_{\a_i}\otimes 1+k_{\a_i}\otimes e_{{\a_i}}\,,\\
\Delta(e_{-\al_i})&=1\otimes e_{-\al_i}+e_{-\al_i}\otimes
k^{-1}_{\a_i}\,,\\
\Delta(k_{\a_i})\,&=\,k_{\a_i}\otimes k_{i}\,.
\end{split}\end{equation}
In the current realization,
$\Uqsln$ is generated by the elements
$\eain{i}{n}$, $\fain{i}{n}$, where $i=1,\ldots,N-1$, $n\in\ZZ$;
$\psi_i^{\pm}[n]$, $i=1,\ldots,N-1$, $n\geq 0$,
$\psi_i^-[0]=\left(\psi_i^+[0]\right)^{-1} $.
They are combined into generating functions
$$
e_{i}(z)\,=\,\sum_{n\in\z}\eain{i}{n}z^{-n}\ ,
\quad f_{i}(z)\,=\,\sum_{n\in\z}\fain{i}{n}z^{-n}\ , \quad
\psi^\pm_i(z)\,=\sum_{n\geq 0}\psi_i^{\pm}[n]z^{\mp n}\,,
$$
which satisfy the following relations:
$$
(z-q^{(\a_i,\al_j)}w)e_{i}(z)e_{j}(w)=
e_{j}(w)e_{i}(z)(q^{(\a_i,\al_j)}z-w)\ ,
$$
$$
(z-q^{-(\a_i,\al_j)}w)f_{i}(z)f_{j}(w)=
f_{j}(w)f_{i}(z)(q^{-(\a_i,\al_j)}z-w)\ ,
$$
$$
\psi_{i}^\pm(z)e_{j}(w)\left(\psi_{i}^\pm(z)\right)^{-1}=
\frac{(q^{(\a_i,\al_j)}z-w)}
{(z- q^{(\a_i,\al_j)}w)}e_{j}(w)\ ,
$$
$$
\psi_{i}^\pm(z)f_{j}(w)\left(\psi_{i}^\pm(z)\right)^{-1}=
\frac{(q^{-(\a_i,\al_j)}z-w)}
{(z- q^{-(\a_i,\al_j)}w)}f_{j}(w)\ ,
$$
$$
\psi_{i}^\mu(z)\psi_{j}^\nu(w)=
\psi_{j}^\nu(w)\psi_{i}^\mu(z)\ ,\quad \mu,\nu=\pm\,,
$$
$$[e_{i}(z),f_{j}(w)]=
\frac{\delta_{ij}\delta(z/w)}{q_{}-q^{-1}_{}}\left(
\psi^+_{i}(z)-\psi^-_{i}(w)\right)\,,
$$
and
\begin{eqnarray*}
\ds \mathop{\rm Sym}\limits_{z_1,z_2}
\sk{e_{i}(z_1)e_{i}(z_2)e_{j}(w) -(q+q^{-1})e_{i}(z_1)
e_{j}(w)e_{i}(z_2)
+e_{j}(w)e_{i}(z_1)e_{i}(z_2) }&=&0\,,\\
\ds \mathop{\rm Sym}\limits_{z_1,z_2}
\sk{f_{i}(z_1)f_{i}(z_2)f_{j}(w) -(q+q^{-1})f_{i}(z_1)
f_{j}(w)f_{i}(z_2)
+f_{j}(w)f_{i}(z_1)f_{i}(z_2) }&=&0\,,
\end{eqnarray*}
where $i-j=\pm 1 $.

The two realizations are related by the formulae:
\begin{align*}
k_{\a_i}&= \psi_i^+[0],\qquad
e_{\a_i}= \eain{i}{0},\qquad e_{-\a_i}= \fain{i}{0},\qquad
i=1,\ldots,N-1,\\
e_{\a_0}&=[e_{1}[0],[e_{2}[0],\ldots,
[e_{{N-2}}[0],e_{{N-1}}[-1]_q]_q\ \ldots\ ]_q\,,\\
e_{-\a_0}&=[\ \ldots\ [f_{{N-1}}[1],f_{{N-2}}[0]]_{q^{-1}},\ldots,
f_{2}[0]]_{q^{-1}},f_{1}[0]]_{q^{-1}}\,,
\end{align*}
where
$[e_{{i}}[k],e_{{j}}[l]]_q=e_{{i}}[k]e_{{j}}[l]-q^{(\a_i,\a_j)}
e_{{j}}[l]e_{{i}}[k]$ \ and \
$[f_{{i}}[k],f_{{j}}[l]]_{q^{-1}} = f_{{j}}[l]f_{{i}}[k] -
q^{-(\a_i,\a_j)}
f_{{i}}[k]f_{{j}}[l]$.

The Drinfeld comultiplication $\Delta^{(D)}$ for the algebra $\Uqsln$
looks as follows.
\begin{align*}
\Delta^{(D)}e_{i}(z)&=e_{i}(z)\otimes 1+\psi_i^-(z)\otimes e_{i} (z)\, ,
\\
\Delta^{(D)}f_{i}(z)&=1\otimes f_{i}(z)+f_{i}(z)\otimes
\psi_i^+(z)\, ,\\
\Delta^{(D)}\psi_i^\pm(z)&=\psi_i^\pm(z)\otimes
\psi_i^\pm(z)\, .
\end{align*}
The quantum affine algebra $\Uqsln$ has two types of Borel subalgebras.
The Borel subalgebras $U_q(\mathfrak{b}_\pm^{\slg})\subset\Uqsln $ are
generated by the Chevalley generators $e_{\a_i}, k^{\pm 1}_{\a_i}$,
$i=0,\ldots,N-1$ and $e_{-\a_i}, k^{\pm 1}_{\a_i}$, $i=0,\ldots,N-1$,
respectively. They contain Hopf coideals $U_q(\mathfrak{n}_\pm^{\slg})\subset
U_q(\mathfrak{b}^{\slg}_\pm)$, generated by the Chevalley generators $e_{\a_i}$,
and
$e_{-\a_i}$, $i=0,\dots,N-1$, respectively.

The Borel subalgebra $U_F^{\slg}\subset \Uqsln$ is generated by the elements
$\fain{i}{n}$, where $i=1,\ldots,N-1$, $n\in\ZZ$ and
$\psi_i^{+}[n]$, $i=1,\ldots,N-1$, $n\geq 0$. The Borel subalgebra
$U_E^{\slg}\subset \Uqsln$ is generated by the elements
$\eain{i}{n}$, where $i=1,\ldots,N-1$, $n\in\ZZ$ and
$\psi_i^{-}[n]$, $i=1,\ldots,N-1$, $n\geq 0$. We are interested in
their intersections,
\begin{equation}
\label{Inter}
U_f^{\slg\,-}=U_F^{\slg}\cap U_q(\mathfrak{n}_-^{\slg})\,,\qquad
U_F^{\slg\,+}=U_F^{\slg}\cap U_q(\mathfrak{b}_+^{\slg})\,.
\end{equation}
According to \cite{EKhP}, these intersections satisfy coideal properties
\begin{equation*}
\Delta^{(D)}(U_F^{\slg\,+})\subset \Uqsln\ot U_F^{\slg\,+},\qquad
\Delta^{(D)}(U_f^{\slg\,-})\subset U_f^{\slg\,-}\ot \Uqsln
\end{equation*}
and the multiplication $m$ in $\Uqsln$ induces an isomorphism of vector spaces
$$m: U_f^{\slg\,-}\ot U_F^{\slg\,+}\to U_F^{\slg}\,.$$
The projection operator $P:U_F^{\slg}\to U_F^{\slg\,+}$ is defined by the rule
\begin{equation}\label{psln}
P(f_-f_+)=\varepsilon(f_-)f_+, \qquad f_-\in U_f^{\slg\,-}\subset \Uqsln,
\qquad f_+\in U_F^{\slg\,+}\subset \Uqsln .
\end{equation}

\subsection{The embedding of $\Uqsln$ to $\Uqgln$}
\label{Thtsect}

Consider the embedding $\Theta:\Uqsln\hookrightarrow\Uqgln$ given by the
following formulae in the current realizations of $\Uqsln$ and $\Uqgln$:
\begin{equation}\label{sltogl}
\begin{split}
\Theta(e_i(z))=(q-q^{-1})^{-1}\ee_i(q^{-i+1}z),\quad&\quad
\Theta(f_i(z))=(q-q^{-1})^{-1}\ff_i(q^{-i+1}z)\ ,\\
\Theta(\psi_i^\pm(z))&= k_i^\pm(q^{-i+1}z)
\left(k_{i+1}^\pm(q^{-i+1}z)\right)^{-1}\!.\quad
\end{split}
\end{equation}
Let $g^+(z),g^-(z)$ be power series with coefficients in $\CC$,
\begin{equation}\label{sl1}
\begin{split}
g^+(z)&=g^+_0+g_1z^{-1}+\ldots+g_nz^{-n}+\ldots,\\
g^-(z)&=g^-_0+g_{-1}z+\ldots+g_{-n}z^{n}+\ldots,
\end{split}
\end{equation}
satisfying the condition
\begin{equation}\label{sl2}
g^+_0g^-_0=1\ .
\end{equation}
A pair $g^\pm(z)$ defines an automorphism $T_{g^+(z),\,g^-(z)}$
of the algebra $\Uqgln$ by the rule
\begin{equation}\label{sl3}
T_{g^+(z),\,g^-(z)}L^+(z)=g^+(z)L^+(z),\qquad
T_{g^+(z),\,g^-(z)}L^-(z)=g^-(z)L^-(z)\, .
\end{equation}
The following facts are well known.
\begin{proposition}\label{sl4}${}$
\begin{itemize}
\item[(i)]
The embedding $\Theta$ is a morphism of Hopf algebras with respect
to any of the comultiplications $\Delta$ or $\Delta^{(D)}$.
\item[(ii)]
The image of $\Theta$ is the subalgebra of invariants of
all automorphisms $T_{g(z^{}),\,\tilde{g}(z)}$ in $\Uqgln$.
\item[(iii)]
The embedding $\Theta$ maps the Borel subalgebras
$U_q(\mathfrak{b}_\pm^{\slg})\subset \Uqsln $ into the corresponding Borel
subalgebras $U_q(\mathfrak{b}_\pm)\subset \Uqgln $, and the current
Borel subalgebra $U_F^{\slg}\subset \Uqsln $ into the current Borel subalgebra
$U_F\subset \Uqgln $.
\end{itemize}
\end{proposition}

\subsection{Projections}
Clearly, the Borel subalgebras in $\Uqsln$ and $\Uqgln$ differ only by Cartan
currents. We have $N$ Cartan currents $k_i^+(z)$, $i=1,\ldots,N$, in
$U_q(\mathfrak{b}_+)\subset \Uqgln $ and in $U_F\subset \Uqgln $, while we have
$N-1$ Cartan currents $\psi_i^+(z)$, $i=1,\ldots,N-1$, in
$U_q(\mathfrak{b}_+^{\slg})\subset \Uqsln $ and in $U_F^{\slg}\subset\Uqsln$,
and $\Theta(\psi_i^+(z))=
k_i^+(q^{-i+1}z) \left(k_{i+1}^+(q^{-i+1}z)\right)^{-1}$.

We can choose modes of one of the current $k_1^+(z)$, as generators of an
abelian subalgebra $A_1$. Then the multiplication in $\Uqgln $ establishes
an isomorphism (of vector spaces) between the Borel subalgebra
$U_q(\mathfrak{b}_+)\subset \Uqgln $ and the tensor product of $A_1$ and the
image of Borel subalgebra $U_q(\mathfrak{b}_+^{\slg})\subset \Uqsln $.
An analogous statement holds for the current Borel subalgebras
$U_F^{\slg}\subset \Uqsln $ and $U_F\subset \Uqgln $.

This observation implies that the algebras
$U_f^-=U'_F\cap U_q(\mathfrak{b}_-)\subset \Uqgln $ and
$U_F^+=U_F\cap U_q(\mathfrak{b}_+)\subset \Uqgln$, see \rf{Intergl},
satisfy the same properties as the analogous subalgebras \rf{Inter} of
$\Uqsln $. Namely,
they are coideals,
\begin{equation*}
\Delta^{(D)}(U_F^+)\subset \Uqgln\ot U_F^+\,,\qquad
\Delta^{(D)}(U_f^-)\subset U_f^-\ot \Uqgln\,,
\end{equation*}
the multiplication $m$ in $\Uqgln$ induces an isomorphism of vector spaces
$$m: U_f^-\ot U_F^+\to U_F\,,$$
and the projection operator
$P:U_F\subset \Uqgln \to U_F^+$, is defined similarly to \rf{psln},
\begin{equation}\label{pgln}
P(f_-f_+)=\varepsilon(f_-)f_+, \qquad f_-\in U_f^-\subset \Uqgln,
\qquad f_+\in U_F^+ \subset \Uqgln.
\end{equation}
Proposition~\ref{sl4} yields that the definitions \rf{psln} and \rf{pgln} are
consistent, that is, for any element $f\in U_F\subset \Uqsln$, we have
\begin{equation}
\label{sl5}
\Theta\bigl(P(f)\bigr)=P\bigl(\Theta(f)\bigr)\,,
\end{equation}
where $P$ in the left hand side is the projection operator \rf{psln} and
$P$ in the right hand side is the projection operator \rf{pgln}.

\subsection{A construction of the weight function}

Let $V$ be a representation of $\Uqgln$ generated by
a singular vector $v$.
Let $I=\{{i_1},\dots,{i_n}\}$ be an ordered $\Pi$-multiset. Set
\begin{equation} \label{WW7a}
w_{V,I}(\{t_i|_{i\in I}\})=
P\left(F_{\i(i_1)}(t_{i_1})\cdots F_{\i(i_n)}(t_{i_n})
\right)\ v.
\end{equation}

\begin{theorem} \label{th22}
A collection of $V$-valued rational functions $w_{V,I}(\{t_i|_{i\in I}\})$,
given by \r{WW7a}, is a $q$-symmetric weight function .
\end{theorem}

\noindent
{\it Proof.}\ \ Due to \rf{sl5} and \rf{sltogl},
$$w_{V,I}(\{t_i|_{i\in I}\})=(q-q^{-1})^{n}
\Theta\left(P\left(f_{\i(i_1)}(\tilde{t}_{i_1})\cdots f_{\i(i_n)}(\tilde{t}_{i_n})
\right)\right)\ v\, ,$$
where $\tilde{t}_{i_1}, \ldots ,\tilde{t}_{i_n}$ are the variables
${t}_{i_1}, \ldots ,{t}_{i_n}$, shifted by some powers of $q$. Therefore,
the collection of functions $w_{V,I}(\{t_i|_{i\in I}\})$ is a $\Uqsln$ weight
function up to a certain shift of variables, and Theorem \ref{th22} is a
particular case of Theorem~4 in \cite{EKhP}. Let us remind that the key
assertion used in the proof of Theorem~4 in \cite{EKhP} is the following
relation for the comultiplications $\Delta$, $\Delta^{(D)}$, and the projection
operator $P$, that holds in $U_q(\widehat\ggg)$ for any simple Lie algebra
$\ggg$: for any element $f\in U_F^{\ggg}$ and any singular vectors $v_1, v_2$,
one has
\begin{equation}\label{main-pro}
\Delta(P(f))\,v_1\otimes v_2 = (P\otimes P)\Delta^{(D)}(f)\,v_1\otimes v_2\,.
\end{equation}
The $q$-symmetry of the weight function $w$ follows from the defining relations
\rf{gln-com}.
\hfill{$\square$}

\medskip
For an ordered $\Pi$-multiset $I=\{i_1,\ldots,i_n\}$ set
\begin{equation}\label{WW7b}
\w^{\rm P}_{V,I}(\{t_i|_{i\in I}\})=
P\left(F_{\i(i_n)}(t_{i_n})\cdots F_{\i(i_1)}(t_{i_1})\right)\,v\,\cdot
\prod_{i\ord j}\beta(t_i,t_j) \prod_{i\in I}
\Lambda_{\i(i)+1}(t_i)\,,
\end{equation}
where $\b(t_i,t_j)$ is defined by \rf{beta}.
Theorem \ref{th22} and Proposition \ref{prop2.3} imply the following statement.

\begin{corollary}\label{corn}
A collection of $V$-valued rational functions
$\w^{\rm P}_{V,I}(\{t_i|_{i\in I}\})$, given
by \rf{WW7b}, is a $q$-symmetric modified weight function.
\end{corollary}

\section{$L$-operators and modified weight functions}
\label{L-op-con}

\subsection{$\Uqgln$ monodromy}

We borrow the construction below from \cite{VT}. We will need only one
$L$-operator, say $L^+(z)$, which we denote as $L(z)$. It generates the Borel
subalgebra $\Uqbp$, see Section \ref{section2.1}.

Let $M$ be a nonnegative integer. Let $L^{(k)}(z)\in \sk{\CC^N}^{\otimes M}$
be the $L$-operator acting as $L(z)$ on $k$-th tensor factor in the product
$\sk{\CC^{N}}^{\otimes M}$ and as the identity operator in all other factors.
Consider a series in $M$ variables
\begin{equation}\label{product}
\ct_{[M]}(u_1,\ldots,u_M)=L^{(1)}(u_1)\cdots L^{(M)}(u_M)
\cdot \ccR^{(M,\ldots,1)}(u_M,\ldots,u_1)
\end{equation}
with coefficients in $\sk{\End(\CC^N)}^{\ot M}\ot\Uqbp$,
where
\begin{equation}\label{Rproduct}
\ccR^{(M,\ldots,1)}(u_M,\ldots,u_1)=
\prodl_{1\leq i<j\leq M} R^{(ji)}(u_j,u_i)
\end{equation}
In the ordered product of $R$-matrices \r{Rproduct} the factor $R^{(ji)}$ is to
the left of the factor $R^{(ml)}$ if $j>m$, or $j=m$ and $i>l$. Note that due
to the Yang-Baxter equation for $L$-operators element \r{product} can be
rewritten as follows,
\begin{equation}\label{product1}
\ct_{[M]}(u_1,\ldots,u_M)=
\ccR^{(M,\ldots,1)}(u_M,\ldots,u_1)\cdot L^{(M)}(u_M)\cdots L^{(1)}(u_1)
\end{equation}

Consider a special
multiset $I_{\bar n}$ labeled by a sequence of non-negative integers
$\bar{n}= \{n_1,n_{2},\ldots,n_{N-1}\}$,
$\bar{n}\in \ZZ_{\geq 0}^{N-1}$, $|\bar n|=n_1+\cdots+n_{N-1}$.
As an ordered set, $I_{\bar n}$ consists of integers $i$, such that
$1\leq i\leq |{\bar n}|$. The colouring map is
\begin{equation}\label{spe-set}
\i(i)=a\in\Pi\quad{\rm for}\quad 1+n_1+\cdots+n_{a-1}\leq i
\leq n_1+\cdots+n_{a}\,.
\end{equation}
Let us change a numeration of the set of variables
$\{t_i|_{ i\in I_{\bar n}}\}$ as
\begin{equation}\label{set1}
\bar{t}_{\bar{n}}\ = \ \{t^a_i\}\ =\ t^1_1,\ldots,t^1_{n_1},\quad t^{2}_1,\ldots,
t^{2}_{n_{2}},\quad \ldots
\quad , \
t^{N-1}_1,\ldots,t^{N-1}_{n_{N-1}}\ .
\end{equation}

Following \cite{VT}, set
\begin{eqnarray}\label{Vit-el}
&\ds\bbb_{\bar n}(\bar t_{\bar n})
=\prod_{a=1}^{N-1}\prod_{1\leq i<j\leq
n_a}\frac{t^a_i-t^a_j}{q^{-1}t^a_i-qt^a_j}\times\\
&\times({\rm tr})^{\ot|\bar n|}\ot {\rm id}) (\ct_{[|\bar
n|]}(t^1_1,\ldots,t^1_{n_1};\ \ldots\ ;
t^{N-1}_1,\ldots,t^{N-1}_{n_{N-1}})\E_{21}^{\ot n_1}\ot\cdots\ot \E_{N,N-1}^{\ot
n_{N-1}}\ot 1).
\nonumber
\end{eqnarray}
Here ${\rm tr}: {\rm End}(\CC^N)\to \CC$ is the standard trace map. The
coefficients of
$\bbb_{\bar n}(\bar t_{\bar n})$ are elements of the Borel subalgebra
$U_q(\mathfrak{b}_+)$.

Let $S_{\bar n}=S_{n_1}\times \cdots
\times S_{n_{N-1}}$ be the direct product of the symmetric groups. The group
$S_{\bar n}$
naturally acts on functions of $t_1^1,\dots t^{N-1}_{n_{N-1}}$ by permutations
of variables with the same superscript, if $\sigma=\sigma^1\times\cdots\times
\sigma^{N-1}\in S_{\bar n}$, then
$$
^\sigma\!\bar t_{\bar n}=(t^1_{\sigma^1(1)},\ldots,t^1_{\sigma^1(n_1)};\ \ldots\ ;
t^{N-1}_{\sigma^{N-1}(1)},\ldots,t^{N-1}_{\sigma^{N-1}(n_{N-1})})\,.
$$
\begin{proposition}\label{sym}
For any $\sigma\in S_{\bar n}$, we have
\begin{equation}\label{symB}
\bbb_{\bar n}(\bar t_{\bar n})=\bbb_{\bar n}(^\sigma\!\bar t_{\bar n})\,.
\end{equation}
\end{proposition}
\noindent
{\it Proof}.
It suffices to prove the claim assuming that $\sigma$ is the product of
a single simple transposition and the identity permutations. In other words,
$\sigma$ permutes just one pair of variables.

Relations \rf{L-op-com}, the Yang-Baxter equation \rf{YBeq} and the inversion
relation \r{unitar} imply that
\begin{equation}\label{Ve1}
\begin{split}
&P^{(i,i+1)}R^{(i,i+1)}(u_i,u_{i+1})\ct_{[M]}(u_1,\ldots,u_i,u_{i+1},\ldots,u_M)=\\
&\quad =
\ct_{[M]}(u_1,\ldots,u_{i+1},u_{i},\ldots,u_M)P^{(i+1,i)}R^{(i+1,i)}(u_{i+1},u_i)\,,
\end{split}
\end{equation}
where $P^{(i,i+1)}$ is the permutation operator,
$P^{(12)}=\sum_{i,j=1}^N \E_{ij}\otimes \E_{ji}$. It is easy to check that
\begin{equation}\label{Ve2}
P^{(12)}\E_{j+1,j}\otimes \E_{j+1,j}=
\E_{j+1,j}\otimes \E_{j+1,j}=\E_{j+1,j}\otimes \E_{j+1,j} P^{(12)}
\end{equation}
and
\begin{equation}\label{Ve3}
\begin{split}
R^{(12)}(u_1,u_2)\E_{j+1,j}\otimes E_{j+1,j}\,=\,\frac{qu_1-q^{-1}u_2}{u_1-u_2}
\ \E_{j+1,j}\otimes \E_{j+1,j}\,={} & \\[6pt]
{}=\,\E_{j+1,j}\otimes \E_{j+1,j} R^{(12)}(u_1,u_2)\, & .
\end{split}
\end{equation}
If $\sigma$ permutes just one pair of variables, then relations
\r{Ve1}\,--\,\r{Ve3} and the cyclic property of the trace yield formula
\rf{symB}. For example, in the simplest nontrivial case
\begin{equation*}
\begin{split}
&\bbb(u_2,u_1)= \frac{u_2-u_1}{q^{-1}u_2-qu_1}\ {\rm
tr}\sk{\ct(u_2,u_1)\E_{j+1,j}\otimes \E_{j+1,j}}=\\
&\quad =
\frac{u_2-u_1}{q^{-1}u_2-qu_1}\ {\rm tr}\sk{
P^{(12)}R^{(12)}(u_1,u_2)
\ct(u_1,u_2)R^{(21)}(u_2,u_1)^{-1}P^{(12)}\E_{j+1,j}\otimes \E_{j+1,j}
}=\\
&\quad =\frac{u_1-u_2}{q^{-1}u_1-qu_2}\ {\rm tr}\sk{\ct(u_1,u_2)\E_{j+1,j}\otimes
\E_{j+1,j}}=
\bbb(u_1,u_2)\,.
\end{split}
\end{equation*}
Proposition~\ref{sym} is proved.
\hfill$\square$

\medskip
Let $V$ be a $\Uqgln$-module generated by a singular vector $v$, cf.~\r{hwv}.
Let $I_{\bar n}$ be a special multiset \r{spe-set}.
Let $\bar t_{\bar n}$ be the set of variables, corresponding to the set
$I_{\bar{n}}$. Set
\begin{equation}\label{swf1}
\w^\bbb_{V,I_{\bar n}}(\{t_i|_{i\in I_{\bar n}}\})=
\bbb_{\bar n}(\bar t_{\bar n})v\,.
\end{equation}

Clearly, special multisets given by the condition \r{spe-set} are
$\Pi$-ordered multisets with an increasing colouring function.
For a $\Pi$-ordered multiset $I$, set
$$n_a=\#\{i\in I\ |\ \i_I(i)=a\}\,,$$
$a=1,\dots N-1$, and $\bar n=(n_1,\dots n_{N-1})$.
Let $\sigma: I\to I_{\bar n}$ be a unique invertible map intertwining the
colouring maps and such that $\sigma(i)\ord\sigma(j)$ iff $\i_I(i)<i_I(j)$,
or $\i_I(i)=i_I(j)$ and $i\ord j$.
Set
\begin{equation}\label{swf2}
\w_{V,I}^\bbb(t_i|_{i\in I})
=\ ^{\sigma,\tilde{\gamma}}
\w_{V,I_{\bar n}}^\bbb(t_i|_{i\in I_{\bar n}})\,.
\end{equation}

\begin{theorem}
A collection of $V$-valued rational functions
$\w^{\bbb}_{V,I}(\{t_i|_{i\in I}\})$, given
by \rf{swf1}, \rf{swf2}, is a $q$-symmetric modified weight function.
\end{theorem}

{\it Proof.}\enspace
The collection $\w^{\bbb}_{V,I}(\{t_i|_{i\in I}\})$ is $q$-symmetric due to
formula \rf{swf2} and Proposition \ref{sym}. Properties \rf{empty} and
\rf{isomorph} of the collection $\w^{\bbb}_{V,I}(\{t_i|_{i\in I}\})$ are
straightforward, and the comultiplication property \rf{weight22} follows from
Theorem~3.6.3 in \cite{VT}.
\hglue0pt\hfill$\square$

\section{A correspondence of the two constructions}
\label{uniq}

The goal of this section is to verify the following statement.
\smallskip

\noindent
{\bf Conjecture. }{\em 
The modified weight functions $\w^{\rm P}$ and $\w^\bbb$, defined respectively
by formulae \r{WW7b} and \r{swf1}, \r{swf2}, coincide.}
\medskip

The Conjecture is equivalent to the following relations.
Let $v$ be a weight singular vector in some $\Uqgln$-module,
cf.~\rf{hwv}. Take a sequence $\bar{n}=\{n_1,\dots,n_{N-1}\}$, and let
$\bar{t}_{\bar{n}}$ be the set of variables \rf{set1}. Then
\begin{equation}
\begin{split}
\label{BBPFF}
\ds\bbb_{\bar n}(\bar t_{\bar n}) \ v=
P\sk{F_{N-1}(t^{N-1}_{n_{N-1}})\cdots
F_{N-1}(t^{N-1}_{1})\cdots
F_{1}(t^{1}_{n_{1}})\cdots
F_{1}(t^{1}_{1})}\,v\,\times{} & \\[4pt]
{}\times\, \prod_{a=1}^{N-1}\left(\prod_{1\leq i<j\leq n_a}
\dfrac{q^{-1}t_i^a-q t_j^a}{t_i^a- t_j^a}\prod_{1\leq i\leq
n_a}\Lambda_{a+1}(t_i^a)\right) &\,\, .
\end{split}
\end{equation}

In this paper we will prove the Conjecture only for the special case,
see~Theorem~\ref{Cor5.1}. Set
$$
R^+(u,v)=\frac{u-v}{qu-q^{-1}v}\ R(u,v)\in{\textrm{End}}(\CC^N\ot\CC^N)\ot
\CC[[{v}/{u}]]
$$
and
$$
R^-(u,v)=\sk{R^+(v,u)^{-1}}^{21}\in{\textrm{End}}(\CC^N\ot\CC^N)\ot
\CC[[{u}/{v}]]\,,
$$
where $R(u,v)$ is defined in \rf{UqglN-R}. Define the {\it evaluation}
representation $\pi_{z}^{(1)}$ of $\Uqgln$ in the coordinate space $\CC^N$
by the rule $\pi_{z}^{(1)}\sk{L^\pm(u)}=R^\pm(u,z)$. We also denote
the representation space of $\pi_{z}^{(1)}$ as $V_{\omega_1}(z)$. The first
coordinate vector
in $\CC^N$ is a weight singular vector. We denote it by $v_{\omega_1}$.

\begin{theorem}\label{Cor5.1}
Let $V$ be a subquotient of the tensor product
$V_{\omega_1}(z_1)\ot \ldots \ot V_{\omega_1}(z_n)$, generated by the singular
vector $v=v_{\omega_1}\ot\ldots \ot v_{\omega_1}$. Then for any ordered
$\Pi$-multiset $I$
\begin{equation}\label{wfe}
\w^{\rm P}_{V}(\{t_i|_{i\in I}\})\ v =
\w^\bbb_{V}(\{t_i|_{i\in I}\})\ v\,.
\end{equation}
\end{theorem}

\noindent
{\it Proof}. Due to the comultiplication properties of the weight functions
it is sufficient to prove Theorem~\ref{Cor5.1} for $n=1$. In this case, the
weight functions generated by the singular vector $v_{\omega_1}$ are nontrivial
only if $|I|<N$ and $\iota_I(I)=\{1,2,\ldots,|I|\}$. Since the weight functions
$\w^{\rm P}$ and $\w^\bbb$ are $q$-symmetric, it is enough to consider the case
$I=\{1,2,\ldots,|I|\}$ with the colouring map $\iota_I(i)=i$ for any $i\in I$.
In the last case, formula \r{wfe} follows from Proposition~\ref{induction}.

\medskip

Let the sequence $\bar n$ be such that $n_1=\cdots=n_k=1$,
$n_{k+1}=\cdots=n_{N-1}=0$. In this case, the set of variables \r{set1} takes
the form:
\begin{equation}\label{sp-set}
\bar t=\{t^1,t^2,\ldots,t^k\}\,,
\end{equation}
where we omit the needless subscript. In addition, write $\bbb_{[k]}(\bar t)$
instead of $\bbb_{\bar n}(\bar t)$. In the described case, formula~\r{BBPFF}
reads as
$$
\bbb_{[k]}(t^1,\ldots,t^{k})\;v\,=\,
P\sk{F_{k}(t^{k})\cdots F_1(t^1) }\prod_{j=1}^{k} \Lambda_{j+1}(t^{j})\ v\ .
$$
and follows from Proposition~\ref{induction}

Say that a vector $v$ is \emph{a singular vector} if
\begin{equation}\label{singE}
E^+_{i+1,i}(z)\cdot v\,=\,0\ ,\qquad i=1,\ldots, N-1\,.
\end{equation}
Equivalently, a vector $v$ is a singular vector if
\begin{equation}\label{singL}
L^+_{ij}(z)\ v=0\,,\qquad 1\le j<i\le N\,.
\end{equation}
Compared with~\r{hwv}, \r{L-op-tr}, we drop here the requirement that
the vector $v$ is an eigenvector of the Cartan currents $k^+_i(z)$,
$i=1,\ldots,N$, and the diagonal entries $L_{ii}(u)$, $i=1,\ldots,N$,
of the $L$-operator. Notice that for a singular vector $v$, we have
$$
k_i(u)\,v\,=\,L_{ii}(u)\,v\,,\qquad i=1,\ldots,N\,,
$$
and $\,L_{ii}(u)L_{jj}(t)\,v=L_{jj}(t)L_{ii}(u)\,v\;$ for any $i,j$.

\begin{proposition}\label{induction}
Let $v$ be a singular vector. Then for any $\,k=1,\ldots,N-1$,
\begin{equation}\label{ind1}
\bbb_{[k]}(t^1,\ldots,t^{k})\;v\,=\,
P\sk{F_{k}(t^{k})\cdots F_1(t^1)}\,L_{k+1,k+1}(t^k)\cdots L_{22}(t^1)\,v\ .
\end{equation}
\end{proposition}

In the rest of the Section we are proving Proposition~\ref{induction}.
The idea of the proof is as follows. Below we will introduce elements
$\bbb_{[l,k]}(t^l,\ldots,t^k)\in U_q(\mathfrak{b}_+)\subset\Uqgln$ such that
$\bbb_{[1,k]}(t^1,\ldots,t^k)=\bbb_{[k]}(t^1,\ldots,t^k)$ and will obtain
relations \r{re-r1} for those elements. We will also consider projections
of partial products of currents, $P\sk{F_{k}(t^{k})\cdots F_l(t^l) }$ and
will obtain relations \r{re-r2} for those projections. The fact that relations
\r{re-r2} and \r{re-r2} are almost the same will allow us to establish
formula \r{ind1}.

\subsection{Recurrence relation for $\bbb_{[l,k]}(t)v$}

For any $l=1,\ldots, k$ introduce an element
$\bbb_{[l,k]}(t^l,\ldots,t^k)\in U_q(\mathfrak{b}_+)\subset \Uqgln$:
\begin{equation}\label{cal1}
\begin{split}
&\bbb_{[l,k]}(t^l,t^{l+1},\ldots,t^k)\,=\,
{\rm tr}_{1,2,\ldots,k-l+1}\left(\ccR^{(k-l+1,\ldots,1)}(t^k,\ldots,t^l)
\right.\times\\
&\qquad\times
\left.L^{(k-l+1)}(t^k)\cdots L^{(2)}(t^{l+1})L^{(1)}(t^l)
\E^{(k-l+1)}_{k+1,k}\cdots
\E^{(2)}_{l+2,l+1}\E^{(1)}_{l+1,l}\right).
\end{split}
\end{equation}
Recall that for any $A\in\End(\CC^N)$ we denote by
$A^{(i)}\in\End(\CC^N\ot\cdots\ot\CC^N)$ the matrix acting as $A$ in
the $i$-th factor of the tensor product $\CC^N\ot\cdots\ot\CC^N$ and
as the identity matrix in all other factors. We also set
\begin{equation}\label{init}
\bbb_{[k+1,k]}(\cdot)\equiv 1 \ .
\end{equation}
It is clear that $\bbb_{[1,k]}(t^1,\ldots,t^k)$ coincides with
$\bbb_{[k]}(t^1,\ldots,t^k)$.

Let $v$ be a singular weight vector.
Recall that $L_{i,i}(t)\,v=k_i(t)\,v=\Lambda_i(t)\,v$.
We will show that the action of the element $\bbb_{[l,k]}(t^l,\ldots,t^k)$ on
the singular vector $v$ can be expressed using a linear combination of products
of the Gauss coordinates $\f_{i,j}(t^{j-1})$ with $l\leq i<j\leq k+1$.
For example, we have $\bbb_{[k,k]}(t^k)=L_{k,k+1}(t^k)$, so that
$\bbb_{[k,k]}(t^k)\,v=\f_{k,k+1}(t^k) L_{k+1,k+1}(t^k)\,v$, and
$$
\bbb_{[k-1,k]}(t^{k-1},t^k)\,=\,L_{k,k+1}(t^k)L_{k-1,k}(t^{k-1})\,+\,
\frac{(q-q^{-1})\,t^k}{t^k-t^{k-1}}\,L_{k-1,k+1}(t^k)L_{k,k}(t^{k-1})
$$
so that
\begin{equation*}
\begin{split}
& \bbb_{[k-1,k]}(t^{k-1},t^k)\,v\,={}\\
& {}=\,\sk{\f_{k,k+1}(t^k)\,\f_{k-1,k}(t^{k-1})+
\frac{(q-q^{-1})\,t^k}{t^k-t^{k-1}}\,\f_{k-1,k+1}(t^k)}
\,L_{k+1,k+1}(t^k)\,L_{k,k}(t^{k-1})\;v\,.
\end{split}
\end{equation*}
To obtain the required presentation in general, we will use the following
statement.
\begin{proposition}\label{tr-rec}
We have
\begin{equation}\label{re-r1}
\begin{split}
\bbb_{[l,k]}(t^l,\ldots,t^k)\,v\,&{}=\,\sum_{m=l+1}^{k+1}
\bbb_{[m,k]}(t^m,\ldots,t^k)\ \f_{l,m}(t^{m-1})\,\times{} \\
&{}\times\,L_{m,m}(t^{m-1})\cdots L_{l+1,l+1}(t^l)\,v\cdot
\prod_{j=l+1}^{m-1}\frac{(q-q^{-1})\,t^j}{t^j-t^{j-1}}\;.
\end{split}
\end{equation}
\end{proposition}
We start the proof of this Proposition from the next lemma.
\begin{lemma}\label{dif-lem}
\begin{equation}\label{one-st}
\begin{split}
\bbb_{[l,k]} & (t^l,\ldots,t^k)\,v\,=\,
\bbb_{[l+1,k]}(t^{l+1},\ldots,t^k)\cdot L_{l,l+1}(t^l)\,v\,+{}\\[3pt]
&{}+\,{\rm tr}_{2,\ldots,k-l+1}
\Bigl(\ccR^{(k-l+1,\ldots,2)}(t^k,\ldots,t^{l+1})
L^{(k-l+1)}(t^k)\cdots L^{(2)}(t^{l+1})\,\times{}\\[-1pt]
&\hphantom{{}+{\rm tr}_{2,\ldots,k-l+1}\Bigl(}{}\times\,
\E^{(k-l+1)}_{k+1,k}\cdots \E^{(3)}_{l+3,l+2}\E^{(2)}_{l+2,l} \Bigr)
\,\frac{(q-q^{-1})\,t^{l+1}}{t^{l+1}-t^l}\ L_{l+1,l+1}(t^l)\,v\ .
\end{split}
\end{equation}
\end{lemma}
\noindent
{\it Proof.}\enspace
To obtain formula~\r{one-st} we calculate the trace over the first copy of
$\CC^N$ in formula~\r{cal1}. Using the Yang-Baxter equation \r{YBeq}, we get
\begin{equation*}
\begin{split}
& \ccR^{(k-l+1,\ldots,2,1)}(t^k,\ldots,t^{l+1},t^l)\,={}\\[4pt]
&{}=\,R^{(2,1)}(t^{l+1},t^l)\cdots R^{(k-l+1,1)}(t^{k-l+1},t^l)
\ccR^{(k-l+1,\ldots,2)}(t^k,\ldots,t^{l+1})\ .
\end{split}
\end{equation*}
Due to relations \r{L-op-tr}, we can write the right hand side
of formula~\r{cal1} applied to the weight singular vector $v$
as a sum of two terms,
\begin{equation}\label{cal2}
\begin{split}
{\rm tr}_{2,\ldots,k-l+1}\biggl(
& {\rm tr}_1 \sk{R^{(21)}(t^{l+1},t^l)\cdots
R^{(k-l+1,1)}(t^k,t^l)\,\E^{(1)}_{ll}}\cdot X\cdot L_{l,l+1}(t^l)\,v\,+{}\\
& \!{}+\,{\rm tr}_1 \sk{R^{(21)}(t^{l+1},t^l)\cdots
R^{(k-l+1,1)}(t^k,t^l)\,\E^{(1)}_{l+1,l}}\cdot X\cdot
L_{l+1,l+1}(t^l)\,v\biggr)\,,
\\[-12pt]
\end{split}
\end{equation}
where
$$
X\,=\,\ccR^{(k-l+1,\ldots,2)}(t^k,\ldots,t^{l+1})
\,L^{(k-l+1)}(t^k)\cdots L^{(2)}(t^{l+1})\,
\E^{(k-l+1)}_{k+1,k}\cdots\E^{(2)}_{l+2,l+1}\,.
$$
Now we calculate the traces $\,{\rm tr}_1$ in formula \r{cal2} taking
into account the matrix structure of the $R$-matrix \r{UqglN-R} and
the multiplication rule: $\E_{ab}\,\E_{cd}=0$ for $b\ne c$ and
$\E_{ab}\,\E_{bc}=\E_{ac}$. As a result, we get
\begin{equation}\label{cal5}
\kern-.4em
{\rm tr}_1\sk{R^{(21)}(t^{l+1},t^l)\cdots R^{(k-l+1,1)}(t^k,t^l)\,
\E^{(1)}_{ll}}\,=\,
\mathbf{1}^{(2)}\,\mathbf{1}^{(3)}\cdots\,\mathbf{1}^{(k-l+1)}+
\sum_{i=2}^{k-l+1}\!\E_{ll}^{(i)}\,Y_i
\end{equation}
and
\begin{equation}\label{cal6}
\begin{split}
\\[-32pt]
{\rm tr}_1\sk{R^{(21)}(t^{l+1},t^l)\cdots R^{(k-l+1,1)}(t^k,t^l)\,
\E^{(1)}_{l+1,l}}\, & {}={}\\[4pt]
{}=\,\frac{(q-q^{-1})\,t^{l+1}}{t^{l+1}-t^l}\ \E_{l+1,l}^{(2)}\,
\mathbf{1}^{(3)}\cdots\,\mathbf{1}^{(k-l+1)}\, &{}+
\sum_{i=3}^{k-l+1}\!\E_{l+1,l}^{(i)}\,Y'_i\,,
\end{split}
\end{equation}
where $Y_i\,, Y'_i$ are some elements of ${\End(\CC^N\ot\cdots\ot\CC^N)}$.
Observe that only the first terms in the right hand sides of formulae~\r{cal5}
and~\r{cal6} contribute nontrivially to the trace $\,{\rm tr}_{2,\ldots,k-l+1}$
in formula~\r{cal2} because ${\rm tr}\,(\E_{ab}\,A\,\E_{cd})=0$
for any $A\in{\rm End}(\CC^N)$ unless $a=d$. To complete the proof of
Lemma~\ref{dif-lem}, notice that
$$
{\rm tr}_{2,\ldots,k-l+1}\,X\,=\,\bbb_{[l+1,k]}(t^{l+1},\ldots,t^k)
$$
and $\;{\rm tr}\,(\E_{l+1,l}\,A\,\E_{l+2,l+1})={\rm tr}\,(A\,\E_{l+2,l})\;$
for any $A\in{\rm End}(\CC^N)$.
\hfill$\square$

\medskip
\noindent
{\it Proof of Proposition~\ref{tr-rec}.}\enspace
To prove this Proposition we use the induction with respect to $N$.

The first term in the right hand side of formula~\r{one-st} is exactly
the term in the right hand side of formula \r{re-r1} for $m=l+1$,
because $L_{l,l+1}(t^l)\,v=F^+_{l,l+1}(t^l)\,L_{l+1,l+1}(t^l)\,v$,
and it suffices to show that
\begin{equation}\label{tr2kl}
\begin{split}
& {\rm tr}_{2,\ldots,k-l+1}
\Bigl(\ccR^{(k-l+1,\ldots,2)}(t^k,\ldots,t^{l+1})
L^{(k-l+1)}(t^k)\cdots L^{(2)}(t^{l+1})\,\times{}\\
&\hphantom{{\rm tr}_{2,\ldots,k-l+1}\Bigl(\kern-1.4em
\ccR^{(k-l+1,\ldots,2)}(t^k,\ldots,t^{l+1})}{}\times\,
\E^{(k-l+1)}_{k+1,k}\cdots \E^{(3)}_{l+3,l+2}\E^{(2)}_{l+2,l}\Bigr)\,v\,={}
\\[3pt]
&{}=\,\sum_{m=l+2}^{k+1}
\bbb_{[m,k]}(t^m,\ldots,t^k)\ \f_{l,m}(t^{m-1})\,
L_{mm}(t^{m-1})\cdots L_{l+2,l+2}(t^{l+1})\,v\cdot
\prod_{j=l+2}^{m-1}\frac{(q-q^{-1})\,t^j}{t^j-t^{j-1}}\;.\kern-1em
\end{split}
\end{equation}

Consider the embedding $\psi:U_q(\widehat{\mathfrak{gl}}_{N-1})\hookrightarrow
\Uqgln$ given by the rule
$$
\psi\sk{L^{[N-1]}_{ij}(t)}=L_{i+\theta(i>l)\,,\,j+\theta(j>l)}(t)\,,
\qquad i,j=1,\ldots,N-1\,,
$$
where $\theta(m>l)=0$ for $m\le l$, and $\theta(m>l)=1$ for $m>l$. Assume that
$U_q(\widehat{\mathfrak{gl}}_{N-1})$ acts by the composition of the embedding
$\psi$ and the action of $\Uqgln$. Then the vector $v$ is singular with respect
to the action of $U_q(\widehat{\mathfrak{gl}}_{N-1})$.
Taking into account the matrix structure of the $R$-matrix \r{UqglN-R},
we can verify that $\psi\sk{\bbb^{[N-1]}_{[m-1,k-1]}(t^m,\ldots,t^k)}=
\bbb_{[m,k]}(t^m,\ldots,t^k)$ for $m>l$, and the left hand side of
formula~\r{tr2kl} coincide with
$\psi\sk{\bbb^{[N-1]}_{[l,k-1]}(t^{l+1},\ldots,t^k)}$. In addition, observe
that
$$
\psi\sk{\bigl(F^+_{l,m-1}\bigr)^{[N-1]}(t)}=F^+_{l,m}(t)\,,\qquad m\ge l+2\,.
$$
As a result, taking formula~\r{re-r1} for $U_q(\widehat{\mathfrak{gl}}_{N-1})$
with parameters $l, k-1, t^{l+1},\dots t^k$, and applying the embedding $\psi$
we obtain formula~\r{tr2kl}. Proposition~\ref{tr-rec} is proved.
\hfill$\square$

\subsection{Composed currents and Gauss coordinates}
\label{MoreDF}

In the next two subsections we will show that the projections of products of
currents, $P\sk{F_{k}(t^{k})\cdots F_l(t^l) }$, satisfy relations \r{re-r2}
which are similar to relations \r{re-r1} for the elements
$\bbb_{[l,k]}(t^l,\ldots,t^k)$. We will use those relations
in subsection~\ref{proof5.1} to prove Proposition~\ref{induction}.

Following \cite{DKh,DKhP,KP}, we will introduce the {\it composed currents\/}
$\ff_{i,j}(t)$ for $i<j$, see \r{res-in}. The composed currents for
an arbitrary quantum affine algebra were defined in \cite{DKh}.
The currents $\ff_{i,j}(t)$ to be used here are images of the composed
currents for the algebra $\Uqsln$ under the embedding
$\Theta:\Uqsln\hookrightarrow\Uqgln$ defined in subsection~\ref{Thtsect}.

The currents $\ff_{i,i+1}(t)$, $i=1,\ldots N-1$, are just the currents
$\ff_i(t)$, cf.~\r{currents}. It follows from formulae \r{DF-iso} and
the definition of the projection $P$, see~\r{pgln}, that
$P\sk{\ff_{i,i+1}(t)}=\f_{i,i+1}(t)$, that is, the projection of the current
$\ff_{i,i+1}(t)$ coincides with the Gauss coordinate $\f_{i,i+1}(t)$
of the corresponding $L$-operator \cite{DF}.
There exists a similar relation between other Gauss coordinates $\f_{i,j}(t)$
and projections of the composed currents $\ff_{i,j}(t)$,
see Proposition~\ref{idenGC}.

According to \cite{DKh}, the composed currents $\ff_{i,j}(t)$ belong to
a suitable completion of the $\Uqgln$ subalgebra generated by modes $F_i[n]$,
$n\in\ZZ$, $i=1,\ldots,N-1$. Elements of the completion are infinite sums of
monomials which are ordered products $\ff_{i_1}[n_1]\cdots \ff_{i_k}[n_k]$ with
$n_1\leq\cdots\leq n_k$. We denote this completion by $\overline U_f$.

The completion $\overline U_f$ determines analyticity properties of products
of currents, see \cite{DKh}. One can show that for $|i-j|>1$, the product
$F_i(t)F_j(w)$ is an expansion of a function analytic at $t\ne 0$, $w\ne 0$.
The situation is more delicate for $j=i,i\pm1$. The products $F_i(t)F_i(w)$ and
$F_i(t)F_{i+1}(w)$ are expansions of analytic functions at $|w|>|q^2 t|$, while
the product $F_i(t)F_{i-1}(w)$ is an expansion of an analytic function at
$|w|>|t|$. Moreover, the only singularity of the corresponding functions
in the whole region $t\ne 0$, $w\ne 0$, are simple poles at the respective
hyperplanes, $w=q^2t$ for $j=i,i+1$, and $w=t$ for $j=i-1$.

The composed currents $F_{i,j}(t)$, $i<j$, are given by the rule
\begin{equation}\label{res-in}
\ff_{i,j}(t)=(q-q^{-1})^{j-i-1}\ff_i(t) \ff_{i+1}(t)\cdots \ff_{j-1}(t)\,.
\end{equation}
For example, $F_{i,i+1}(t)=F_i(t)$, and
$F_{i,i+2}(t)=(q-q^{-1})F_i(t)F_{i+1}(t)$. The last product is well-defined
according to the analyticity properties of the product $F_i(t)F_{i+1}(w)$,
described above. In a similar way, one can show inductively that the product
in the right hand side of \r{res-in} makes sense for any $i<j$.

Products of the composed currents have the following analyticity properties.
For any $i<r<s<j$, the products $F_{i,r}(t)F_{s,j}(w)$ and $F_{j,s}(t)F_{i,r}(w)$
are expansions of functions analytic at $t\ne 0$, $w\ne 0$. For any $i<s<j$,
the product $F_{i,s}(t)F_{s,j}(w)$ is an expansion of an analytic function
at $|w|>|q^2 t|$, and the product $F_{s,j}(t)F_{i,s}(w)$ is an expansion of
an analytic function at $|w|>|t|$. Moreover, the only singularity of the
corresponding functions in the whole region $t\ne 0$, $w\ne 0$, are simple
poles at the respective hyperplanes, $w=q^2t$ for $F_{i,s}(t)F_{s,j}(w)$,
and $w=t$ for $F_{s,j}(t)F_{i,s}(w)$.

The composed currents obey commutation relations
\begin{equation}\label{5.61}
(q^{-1}w-qt)\ff_{i,s}(w)\ff_{s,j}(t)=(w-t)\ff_{s,j}(t)\ff_{i,s}(w)\,,
\end{equation}
for any $i<s<j$, and
\begin{equation}\label{5.62}
\ff_{i,r}(w)\ff_{s,j}(t)=\ff_{s,j}(t)\ff_{i,r}(w)\,,
\end{equation}
for any $i<r<s<j$, which can be observed from the basic relations \r{gln-com}
and formula~\r{res-in}.
In addition, the residue formula
\begin{equation}\label{rec-f1}
\ff_{i,j}(t)\, =\,
-\,\mathop{\rm res}\limits_{w=t}\ff_{s,j}(t)\ff_{i,s}(w)\,\frac{dw}{w}
\end{equation}
holds for any $s=i+1,\ldots, j-1$. Since the total sum of residues of an
analytic functions equals zero, taking into account commutation relations
\r{5.61}, we also get
\begin{equation}\label{rec-f111}
\ff_{i,j}(t)\,=\,
\mathop{\rm res}\limits_{w=0}\sk{ \ff_{s,j}(t) \ff_{i,s}(w)\ \frac{dw}{w}}\,+
\,\mathop{\rm res}\limits_{w=\infty}\sk{ \frac{q^{-1}w-qt}{w-t}
\;\ff_{i,s}(w) \ff_{s,j}(t)\ \frac{dw}{w}}\,.
\end{equation}

Set $S_A(B)=BA-qAB$.
Projections of composed currents can be defined using $q$-commutators
with zero modes of the currents $F_i(t)$, $i=1,\ldots,N-1$.
We will call operators $S_{F_i[0]}$ {\it the screening operators}.
\begin{proposition}\label{com-cur2}
We have
\begin{equation}\label{pro-in}
P\sk{\ff_{i,j}(t)}\,=\,
S_{\ff_i[0]}\bigl(P(\ff_{i+1,j}(t))\bigr)\,,\qquad i<j-1\,.
\end{equation}
\end{proposition}
\noindent
{\it Proof.}\enspace
Calculating the residues in the right hand side of formula \r{rec-f111}
for $s=i+1$ and using the fact $F_{i,i+1}(t)=F_i(t)$, we obtain
\begin{equation}\label{comf-m}
F_{i,j}(t)\,=\,
F_{i+1,j}(t)F_i[0]-qF_i[0]F_{i+1,j}(t)+
(q-q^{-1})\sum_{k\leq 0}F_i[k]\,F_{i+1,j}(t)\,t^{-k}\,.
\end{equation}
Now we apply the projection $P$, see~\r{pgln}, to both sides of this relation.
The modes $F_i[k]$ with $k\leq0$ belong to $U^-_f$. Hence, due to
formulae~\r{pgln}, the projection $P$ kills the semi-infinite sum in
the right hand side of \r{comf-m}, and we get
\begin{equation}\label{pro1}
\begin{split}
P\sk{F_{i,j}(t)}\,&{}=\,
P\sk{F_{i+1,j}(t)F_i[0]-q\,F_i[0]F_{i+1,j}(t)}\,={}\\[4pt]
&{}=\,P\sk{S_{F_i[0]}\sk{F_{i+1,j}(t)}}\,=\,
S_{F_i[0]}\sk{P\sk{F_{i+1,j}(t)}}
\end{split}
\end{equation}
To get the last equality we use the fact proved in \cite{DKhP,KP} that
the projection $P$ commutes with the screening operators $S_{F_i[0]}$,
$P\bigl(S_{F_i[0]}(F)\bigr)=S_{F_i[0]}\bigl(P(F)\bigr)$
for any $F\in U'_F$.
\hfill$\square$

\begin{proposition}\label{idenGC}
We have
\begin{equation}\label{ide-pr}
P\sk{F_{i,j}(t)}\,=\,(q-q^{-1})^{j-i-1}\f_{i,j}(t)\,,\qquad i<j-1\,.
\end{equation}
\end{proposition}
\noindent{\it Proof.}\enspace
The claim follows by induction with respect to $j-i$ from formula
$P\bigl(F_{i,i+1}(u)\bigr)=F^+_{i,i+1}(u)$, Proposition~\ref{com-cur2}
and Lemma~\ref{5.6} proved below.
\hfill$\square$

\begin{lemma}\label{5.6}
We have
\begin{equation}\label{cor1}
(q-q^{-1})\f_{i,j}(t)\,=\,S_{F_i[0]}\sk{\f_{i+1,j}(t)}\,,\qquad i<j-1\,.
\end{equation}
\end{lemma}
\noindent
{\it Proof.}\enspace
It follows from relation \r{L-op} that
\begin{equation}\label{LF}
L^+_{i,j}(t)=
\f_{i,j}(t)k^+_j(t)+\sum_{m=j+1}^N \f_{i,m}(t)k^+_m(t)\e_{m,j}(t)\,.
\end{equation}
Since $S_{F_i[0]}\sk{\f_{i+1,j}(t)}=\f_{i+1,j}(t)F_i[0]-qF_i[0] \f_{i+1,j}(t)$,
and taking into account the commutativity
$$
[F_i[0],k^+_{j}(t)]=0\qquad{\rm and}\qquad
[F_i[0],k^+_m(t)\e_{m,j}(t)]=0
$$
for $i=1,\ldots,j-2$ and $m=j+1,\ldots,N$, we observe that relation~\r{cor1}
results from formula~\r{LF} and the equality
\begin{equation}\label{cor3}
(q-q^{-1})L^+_{i,j}(t)=L^+_{i+1,j}(t)F_i[0]-q F_i[0] L^+_{i+1,j}(t)
\end{equation}
by induction with respect to $j$ starting from $j=N$. On the other hand,
the second line in~\r{L-op-com} at $w=0$ gives
$$
L^+_{i+1,j}(t)L^-_{i,i+1}[0]+(q-q^{-1})L^+_{i,j}(t)L^-_{i+1,i+1}[0]\,=\,
L^-_{i,i+1}[0]L^+_{i+1,j}(t)\,,
$$
which yields formula~\r{cor3}, if we keep in mind the relations
\begin{equation*}
\begin{split}
L^-_{i,i+1}(0)=L^-_{i,i+1} & [0]=-F_i[0]k_{i+1}^{-1}\,,\qquad
L^-_{i,i}(0)=L^-_{i,i}[0]=k_i^{-1}\,,\\[4pt]
& k^{-1}_{i+1}L^+_{i+1,j}(t)k_{i+1}=qL^+_{i+1,j}(t)\,,
\end{split}
\end{equation*}
also following from~\r{L-op-com}, \r{L-op}.
\hfill$\square$

\subsection{Calculation of the projections}

The following proposition is a counterpart of Proposition~\r{tr-rec}.

\begin{proposition}\label{rec-pro}
For any $k>l$, we have
\begin{equation}\label{re-r2}
\begin{split}
&P\sk{F_{k,k+1}(t^{k})\cdots F_{l,l+1}(t^l) }=\\
&\qquad =\sum _{m=l+1}^{k+1}
P\sk{F_{k,k+1}(t^{k})\cdots F_{m,m+1}(t^m) }
P\sk{F_{l,m}(t^{m-1})}\prod_{j=l+1}^{m-1}\frac{t^{j}}
{t^{j}-t^{j-1}}\,.
\end{split}
\end{equation}
\end{proposition}
\noindent
{\it Proof.\/}\enspace
The claim follows from Lemma~\ref{steps} proved below.
\hfill$\square$.

\begin{lemma}\label{steps}
For any $j=l+1,\ldots,k$, we have
\begin{equation}\label{alt1}
\begin{split}
&\kern-1em
P\sk{F_{k,k+1}(t^k)F_{k-1,k}(t^{k-1})\cdots
F_{j+1,j+2}(t^{j+1})F_{j,j+1}(t^{j})F_{l,j}(t^{j-1})}\,={}\\[6pt]
&{}=\,P\sk{F_{k,k+1}(t^k)F_{k-1,k}(t^{k-1})\cdots
F_{j+1,j+2}(t^{j+1})F_{j,j+1}(t^{j})} P\sk{F_{l,j}(t^{j-1})}\,+{}\\
&\qquad\ {}+\,\frac{t^j}{t^j-t^{j-1}}\ P\sk{F_{k,k+1}(t^k)
F_{k-1,k}(t^{k-1})\cdots F_{j+1,j+2}(t^{j+1})F_{l,j+1}(t^{j})}\,.
\end{split}
\end{equation}
\end{lemma}
{\it Proof.\/}\enspace
We will prove the lemma by induction with respect to $l$,
decreasing $l$ from $j-1$ to $1$.

We will use the following properties of the projection $P$ coming from
the definition~\r{pgln}:
\begin{equation}\label{alt3}
P(F^-_{s,s+1}(t)\cdot X)\,=\,0\,,
\end{equation}
for any $s=1,\ldots N-1$, and $X\in U_f\,$, and
\begin{equation}\label{PXX}
P(X_1P(X_2))=P(X_1)P(X_2)\,,
\end{equation}
for any $X_1,X_2\in U_f\,$.

Consider the case $l=j-1$. To calculate the left hand side of
formula \r{alt1} we split the current $F_{l,j}(t^{j-1})=F_{j-1,j}(t^{j-1})$,
$$
F_{j-1,j}(t^{j-1})\,=\,F^+_{j-1,j}(t^{j-1})-F^-_{j-1,j}(t^{j-1})\,.
$$
Since $F^+_{j-1,j}(t^{j-1})=P\bigl(F_{j-1,j}(t^{j-1})\bigr)$,
the first term here produces the first term in the right hand side of
formula~\r{alt1}. For the second term, we move the negative half-current
$F^-_{j-1,j}(t^{j-1})$ to the left using the relation
\begin{equation}\label{cpr2}
\begin{split}
F_{j,j+1}(t^j)F^-_{j-1,j}(t^{j-1}) \, &
{}=\,\frac{qt^j-q^{-1}t^{j-1}}{t^j-t^{j-1}}\;F^-_{j-1,j}(t^{j-1})F_{j,j+1}(t^j)\,-{}\\[4pt]
&{}-\,\frac{(q-q^{-1})\,t^j}{t^j-t^{j-1}}\;F^-_{j-1,j}(t^j)F_{j,j+1}(t^j)
+\frac{t^j}{t^j-t^{j-1}}\;F_{j-1,j+1}(t^j)\,,
\end{split}
\end{equation}
which is a consequence of formulae~\r{5.61}, \r{rec-f1} and the analyticity
properties of the products of currents, and the fact that the currents
$F_{s,s+1}(t^s)$ for $s>j$ commute with $F^-_{j-1,j}(t^{j-1})\,$.
Due to relation~\r{alt3}, only the third term in the right hand side
of~\r{cpr2} contributes nontrivially to the projection, and we obtain
the second term in the right hand side of formula~\r{alt1} for $l=j-1$,
\begin{equation}\label{cpr3}
\begin{split}
&P\sk{F_{k,k+1}(t^{k})\cdots F_{l+1,l+2}(t^{l+1}) }\cdot P\sk{F_{l,l+1}(t^l)}
\,+{}\\
&{}+\,P\sk{F_{k,k+1}(t^{k})\cdots F_{l+2,l+3}(t^{l+2}) F_{l,l+2}(t^{l+1})}\
\frac{t^{l+1}}{t^{l+1}-t^l}\,.
\end{split}
\end{equation}

Assume now that $l\le j-2$. Formula~\r{comf-m} gives that
\begin{equation}\label{eq1}
F_{l,s}(t)\,=\,
S_{F_l[0]}\sk{F_{l+1,s}(t)} - (q-q^{-1})F^-_{l,l+1}(t)\,F_{l+1,s}(t)\,,
\end{equation}
We replace $F_{l,j}(t^{j-1})$ and $F_{l,j+1}(t^j)$ in~\r{alt1} by the right
hand side of formula~\r{eq1} for $s=j,j+1$, respectively. Since the currents
$F_{s,s+1}(t^s)$ for $s>l+1$ commute with $F^-_{l,l+1}(t)\,$, the contributions
of the second term in~\r{eq1} vanish due to relation~\r{alt3}.
For the first term, we use the fact that
$P\bigl(S_{F_l[0]}(F)\bigr)=S_{F_l[0]}\bigl(P(F)\bigr)$ for any $F\in U'_F$,
see~\cite{DKhP,KP}, and the commutativity of the currents $F_{s,s+1}(t^s)$ for
$s>l+1$ with $F_l[0]\,$. As a result, we get that formula~\r{alt1} is
equivalent to
\begin{equation*}
\begin{split}
&\kern-1em S_{F_l[0]}\bigl(P\sk{F_{k,k+1}(t^k)F_{k-1,k}(t^{k-1})\cdots
F_{j+1,j+2}(t^{j+1})F_{j,j+1}(t^{j})F_{l+1,j}(t^{j-1})}\bigr)\,={}\\[6pt]
&{}=\,S_{F_l[0]}\bigl(P\sk{F_{k,k+1}(t^k)F_{k-1,k}(t^{k-1})\cdots
F_{j+1,j+2}(t^{j+1})F_{j,j+1}(t^{j})}P\sk{F_{l+1,j}(t^{j-1})}\bigr)\,+{}\\
&\qquad\quad{}+\,\frac{t^j}{t^j-t^{j-1}}
\ S_{F_l[0]}\bigl(P\sk{F_{k,k+1}(t^k)F_{k-1,k}(t^{k-1})\cdots
F_{j+1,j+2}(t^{j+1})F_{l+1,j+1}(t^{j})}\bigr)\,.
\end{split}
\end{equation*}
The last equality is obtained by application of the screening operator
$S_{F_l[0]}$ to formula~\r{alt1} with $l$ replaced by $l+1$, and is true
by the induction assumption.
\hfill$\square$

\subsection{Proof of Proposition~\ref{induction}}
\label{proof5.1}

For each $l=1,\ldots N-2$, we consider the embedding
$\psi_l:U_q(\widehat{\mathfrak{gl}}_{N-l})\hookrightarrow\Uqgln$,
defined by the rule
$$
\psi_l\sk{L^{[N-l]}_{ij}(t)}=L_{i+l\,,\,j+l}(t)\,,
\qquad i,j=1,\ldots,N-l\,.
$$
Taking into account the matrix structure of the $R$-matrix \r{UqglN-R},
one can verify that for $l<m<k$,
$$
\psi_l\sk{\bbb^{[N-l]}_{[m-l,k-l]}(t^m,\ldots,t^k)}\,=\,
\bbb_{[m,k]}(t^m,\ldots,t^k)\,.
$$
In addition, using formula~\r{L-op} one can check that for $l<i$,
\ $\psi_l\sk{F^{[N-l]}_{i-l}(t)}\,=\,F_i(t)$.
Besides this, the embedding $\psi_l$ is consistent with the projections
$P^{[N-l]}$ and $P$.

We prove Proposition~\r{induction} by induction with respect to $N$.
We replace the expressions in both sides of formula~\r{ind1} by the right hand
sides of formulae~\r{re-r1} and~\r{re-r2} with $l=1$, respectively, and compare
the results term by term. The terms for $m=k+1$ are manifestly the same, taking
into account formula~\r{ide-pr}. For $m=2,\ldots k$, the equality of the
corresponding terms is equivalent to
\begin{equation}
\label{psiBPF}
\begin{split}
\psi_{m-1} & \sk{\bbb^{[N-m+1]}_{[k-m+1]}(t^m,\ldots,t^k)}
\f_{1,m}(t^{m-1})\,L_{mm}(t^{m-1})\cdots L_{22}(t^1)\,v\,={}\\[4pt]
{}=\,{} & \psi_{m-1}\sk{P^{[N-m+1]}\bigl(F^{[N-m+1]}_{k-m+1}(t^k)\cdots
F^{[N-m+1]}_1(t^m)\bigr)}\,\times{}\\[4pt]
& \hphantom{\psi_{m-1}\sk{P^{[N-m+1]}\,F^{[N-m+1]}}}
\times\,\f_{1,m}(t^{m-1})\,L_{k+1,k+1}(t^k)\cdots L_{22}(t^1)\,v\,.
\end{split}
\end{equation}

It follows from commutation relations~\r{L-op-com}, \r{gln-com}, that if $v$
is a singular vector with respect to the action of $\Uqgln$, then the vector
$$
v_{m-1}\,=\,\f_{1,m}(t^{m-1})\,L_{mm}(t^{m-1})\cdots L_{22}(t^1)\,v
$$
is a singular vector with respect to the action of
$U_q(\widehat{\mathfrak{gl}}_{N-m+1})$ induced by the embedding
$\psi_{m-1}:U_q(\widehat{\mathfrak{gl}}_{N-m+1})\hookrightarrow\Uqgln\,$,
and
$$
\f_{1,m}(t^{m-1})\,L_{k+1,k+1}(t^k)\cdots L_{22}(t^1)\,v\,=\,
L_{k+1,k+1}(t^k)\cdots L_{m+1,m+1}(t^m)\,v_{m-1}\,.
$$
Hence, formula~\r{psiBPF} takes the form
\begin{equation*}
\begin{split}
\psi_{m-1} & \sk{\bbb^{[N-m+1]}_{[k-m+1]}(t^m,\ldots,t^k)}\,v_{m-1}\,={}\\[4pt]
{}=\,{} & \psi_{m-1}\Bigl(P^{[N-m+1]}\bigl(F^{[N-m+1]}_{k-m+1}(t^k)\cdots
F^{[N-m+1]}_1(t^m)\bigr)\times{}\\[4pt]
& \hphantom{\psi_{m-1}\,F^{[N-m+1]}}
\times\,L_{k-m+2,k-m+2}(t^k)\cdots L_{22}(t^m)\Bigr)\,v_{m-1}\,,
\end{split}
\end{equation*}
which follows from the induction assumption.
\hfill$\square$

\section*{Acknowledgement}
The work of S.\,Khoroshkin and S.\,Pakuliak was supported in part by grants
INTAS-OPEN-03-51-3350, RFBR grant 04-01-00642 and RFBR grant NSh-8065.2006.2
to support scientific schools. The work of V.\,Tarasov was supported in part
by RFBR grant 05-01-00922. This work was partially done when the first two
authors visited the Max Plank Institut f\"ur Mathematik in Bonn. They thank
the MPIM for the hospitality and stimulating scientific atmosphere.


\begin{thebibliography}{DKhP}

\bibitem[D]{D88} Drinfeld, V. New realization of Yangians and quantum
affine algebras. {\it Sov. Math. Dokl.} {\bf 36} (1988) 212--216.

\bibitem[DF]{DF} Ding, J., Frenkel, I.B. Isomorphism of two realizations of
quantum affine algebra $\Uqgln$. Comm. Math. Phys. \textbf{156} (1993),
277--300.

\bibitem[DKP1]{DKhP1} Ding, J., Khoroshkin, S., Pakuliak, S. Integral
presentations for the universal $R$-matrix. {\it Lett. Math. Phys.} {\bf 53}
(2000), no. 2, 121--141.

\bibitem[DK]{DKh} Ding, J., Khoroshkin, S. Weyl group extension of quantized
current algebras. {\it Transformation Groups}. {\bf 5} (2000), 35--59.

\bibitem[DKP]{DKhP} Ding, J., Khoroshkin, S., Pakuliak, S. Factorization of
the universal $R$-matrix for $U_q(\widehat{sl}_2)$ {\it Theor. and Math. Phys.}
{\bf 124:2} (2000), 1007-1036.

\bibitem[E]{E} B. Enriquez, On correlation functions of Drinfeld
currents and shuffle algebras, Transform. Groups 5 (2000), n.2,
111-120.

\bibitem[EKP]{EKhP} Enriquez, B., Khoroshkin, S., Pakuliak, S. Weight
functions and Drinfeld currents. Preprint ITEP-TH-40/05, {\tt math.QA/0610398}.

\bibitem[ER]{ER} Enriquez, B., Rubtsov, V. Quasi-Hopf algebras associated with
$\mathfrak{sl}_2$ and complex curves. {\it Israel J. Math} {\bf 112} (1999)
61--108.

\bibitem[KP]{KP} Khoroshkin, S., Pakuliak, S. Weight function for
$U_q(\widehat{\mathfrak{sl}}_3)$ {\it Theor. and
Math. Phys.}, {\bf 145} (2005), no.~1, 1373--1399, {\tt math.QA/0610433}.

\bibitem[KT]{KT} Khoroshkin, S., Tolstoy, V. Twisting of quantum
(super)algebras. Connection of Drinfeld's and Cartan-Weyl realizations
for quantum affine algebras. {\it MPI Preprint MPI/94-23},
{\tt hep-th/9404036}.

\bibitem[KR]{KR83} Kulish, P., Reshetikhin, N. Diagonalization of $GL(N)$
invariant transfer matrices and quantum $N$-wave system (Lee model)
{\sl J.Phys. A: Math. Gen.} {\bf 16} (1983) L591--L596

\bibitem[R]{R} Reshetikhin, N. Jackson-type integrals, Bethe vectors, and
solutions to a difference analogue of the Knizhnik-Zamolodchikov system.
Lett. Math. Phys. \textbf{26} (1992), 153--165.

\bibitem[RS]{RS} Reshetikhin, N., Semenov-Tian-Shansky, M. Central extentions
of quantum current groups. Lett. Math. Phys. \textbf{19} (1990), 133--142.

\bibitem[S]{S} Smirnov, F. Form factors in completely integrable models of
quantum field theory, Adv. Series in Math. Phys., vol. 14, World Scientific,
Singapore, 1992.

\bibitem[TV1]{VT} Tarasov, V., Varchenko, A. Jackson integrals for the
solutions to Knizhnik-Zamolodchikov equation, {\it Algebra and Analysis} {\bf
2} (1995) no.2, 275--313.

\bibitem[TV2]{VT1} Tarasov, V., Varchenko, A.
Geometry of $q$-hypergeometric functions, quantum affine algebras
and elliptic quantum groups, {\sl Ast\'erisque} {\bf 246} (1997), 1--135.

\end{thebibliography}
\end{document}